\documentclass[11 pt,Rectoverso]{article}  
\def \cw{\check{w}}
\def \Fbin{{\cal F}_{\sf bin}}
\def \Fter{{\cal F}_{\sf ter}}

\def\di{/\!\!/}

\usepackage[latin1]{inputenc}
\usepackage{graphicx,psfrag,amsfonts,bbm,color}
\usepackage{amsmath,indentfirst,qsymbols}

\def \app#1#2#3#4#5{\begin{array}{rccl} #1:&#2&\longrightarrow&#3\\ &#4&\longmapsto&#5\end{array}}

\def \Geo{{\sf Geometric}}

\def \1{\mathbbm{1}}

\def \bin{{\sf bin}}
\def \ter{{\sf ter}}
\def \thin{{\sf thin}}
\def \Tbin{{\cal T}^{\sf bin}}
\def \Tter{{\cal T}^{\sf ter}}
\def \Tterb{{\cal T}^{\sf ter \bullet}}
\def \Tters{{\cal T}^{\sf ter \star}}
\def \be{\begin{eqnarray*}}
\def \ee{\end{eqnarray*}}
\def \ben{\begin{eqnarray}}
\def \een{\end{eqnarray}}
\def \bq{\begin{equation}}
\def \eq{\end{equation}}

\def \build#1#2#3{\mathrel{\mathop{\kern 0pt#1}\limits_{#2}^{#3}}}
\def \cro#1{\llbracket#1\rrbracket}

\def \se{{\sf e}}

\def \floor#1{\lfloor#1\rfloor}

\def \eref#1{(\ref{#1})}
\def \sous#1#2{\mathrel{\mathop{\kern 0pt#1}\limits_{#2}}}
\def \sur#1#2{\mathrel{\mathop{\kern 0pt#1}\limits^{#2}}}

\def \dd{\xrightarrow[n]{(d)}}
\def \as{\xrightarrow[n]{(a.s.)}}
\def \proba{\xrightarrow[n]{proba.}}

\def \proof {\noindent\bf Proof. \rm}
\def \prooff#1{\noindent\bf Proof of {#1}. \rm}

\def \sa{{\sf a}}
\def \captionn#1{\begin{center}\begin{minipage}{14cm}\sf\caption{\small #1}\end{minipage}\end{center}}

\def \T{{\cal T}}

\def \tend{\longrightarrow}

\def \B{\Big}

\def \l{\left}
\def \r{\right}

\def \B{B}

\def \N{\mathbb{N}}

\def \bt{{\bf t}}

\font\dsrom=dsrom10 scaled 1400
\newqsymbol{`I}{\textrm{\dsrom{1}}}
\newqsymbol{`P}{\mathbb{P}}
\newqsymbol{`R}{\mathbb{R}}
\newqsymbol{`Z}{\mathbb{Z}}
\newqsymbol{`E}{\mathbb{E}}
\newqsymbol{`e}{\varepsilon}
\newqsymbol{`a}{\alpha}
\newqsymbol{`t}{\tau}
\newqsymbol{`T}{\mathbb{T}}
\newqsymbol{`o}{\omega}
\newqsymbol{`O}{\Omega}
\newqsymbol{`M}{\overrightarrow{\cal M}}
\newqsymbol{`m}{{M_n^{\bullet}}}
\newqsymbol{`Q}{\mathbb{Q}}
\newqsymbol{`q}{{{\cal Q}_n^{\bullet}}}
\newqsymbol{`W}{{\cal W}^+}
\newqsymbol{`w}{{\cal W}}
\newqsymbol{`d}{{\mathcal D}}
\newqsymbol{`N}{\mathbb{N}}
\newqsymbol{`Z}{\mathbb{Z}}
\def \ttt{{\bf \triangle}}
\def \3#1{{\bf \triangle}_{#1}}
\def \4#1{{\bf \square}_{#1}}
\newqsymbol{`4}{{\bf \square}}
\newqsymbol{`3}{{\bf \triangle}}
\def \uter{\mathbb{U}^{\ter{}}}
\def \ut{\mathbb{U}^{\3{}}}
\def \uq{\mathbb{U}^{\4{}}}

\def \qt{\mathbb{Q}^{\3{}}}
\def \qbin{\mathbb{Q}^{\bin}}
\def \qter{\mathbb{Q}^{\ter}}
\def \qq{\mathbb{Q}^{\4{}}}

\graphicspath{{../DESSINS/}}
\DeclareMathOperator{\Prof}{Prof}

\DeclareMathOperator{\Dir}{Dir}
\DeclareMathOperator{\type}{type}
\DeclareMathOperator{\Dec}{Dec}

\begin{document}
\renewcommand{\baselinestretch}{1.2}
\newtheorem{exe}{Exercise}
\newtheorem{fig}{\hspace{1cm} Figure}
\newtheorem{lem}{Lemma}
\newtheorem{conj}{Conjecture}
\newtheorem{defi}{Definition}
\newtheorem{pro}[lem]{Proposition}
\newtheorem{theo}[lem]{Theorem}
\newtheorem{cor}[lem]{Corollary}
\newtheorem{remi}{Remark\rm}{\rm}

\newtheorem{com}{Comments\rm}{\rm}
\newenvironment{rem}%
{\begin{center}\begin{minipage}{16cm}\begin{remi}}%
{\end{remi}\end{minipage}\end{center}}
\newcounter{aaa}

\newtheorem{note}{Note \rm}{\rm}

\begin{center}
{\LARGE\bf Some families of increasing planar maps}
\end{center}

\[\begin{array}{c}
\begin{array}{lcl}
\textrm{\Large Marie Albenque}&~~&\textrm{\Large Jean-Fran\c{c}ois Marckert}\\
\textrm{LIAFA, CNRS UMR 7089}&& \textrm{CNRS, LaBRI, UMR 5800}\\
\textrm{Université Paris Diderot - Paris 7}&& \textrm{Universit\'e Bordeaux 1}\\
\textrm{75205 Paris Cedex 13}&&  \textrm{351 cours de la Libération}\\
   &&\textrm{33405 Talence cedex}
\end{array} 

\end{array}\]

\begin{abstract} Stack-triangulations appear as natural objects when one wants to define some increasing families of triangulations by successive additions of faces.
We investigate the asymptotic behavior of rooted stack-triangulations with $2n$ faces under two different distributions. We show that the uniform distribution on this set of maps converges, for a topology of local convergence, to a distribution on the set of infinite maps. In the other hand, we show that rescaled by $n^{1/2}$, they converge for the Gromov-Hausdorff topology on metric spaces to the continuum random tree introduced by Aldous. Under a distribution induced by a natural random construction, the distance between random points rescaled by $(6/11)\log n$ converge to 1 in probability. \par
We obtain similar asymptotic results for a family of increasing quadrangulations.
\end{abstract}

\section{Introduction}

Consider a rooted triangulation of the plane. Choose a finite triangular face $ABC$ and add inside a new vertex $O$ and the three edges $AO$, $BO$ and $CO$. Starting at time 1 from a single rooted triangle, after $k$ such evolutions, a triangulation with $2k+2$ faces is obtained. The set of triangulations  $\3{2k}$  with $2k$ faces that can be reached by this growing procedure is not the set of all rooted triangulations with $2k$ faces. The set $\3{2k}$ -- called the set of stack-triangulations with $2k$ faces -- can be  naturally endowed with two very different probability distributions:
\begin{itemize}
\item[-] the first one, very natural for the combinatorial point of view, is the uniform distribution $\ut_{2k}$, 
\item[-] the second probability $\qt_{2k}$ maybe more realistic following the description given above, is the probability induced by the construction when the faces where the insertion of edges are done are chosen uniformly among the existing finite faces.
\end{itemize}

 \begin{figure}[htbp]
\centerline{\includegraphics[height= 2.5cm]{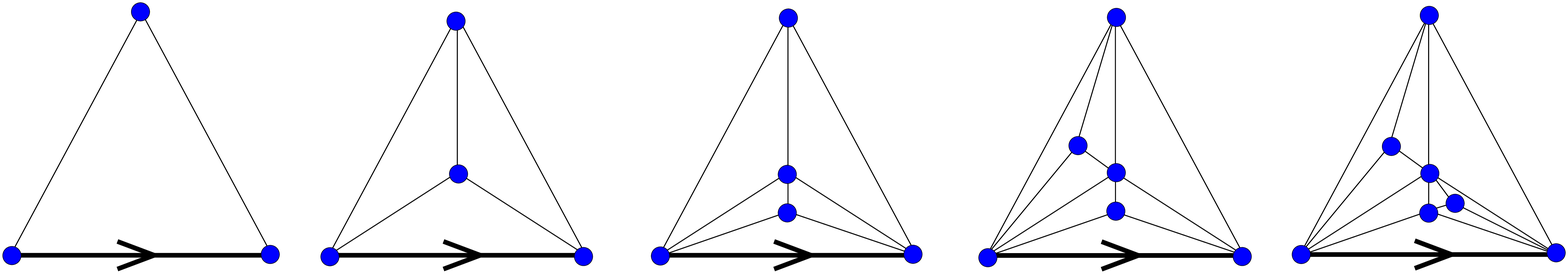}}
\captionn{\label{fig2}Iterative construction of a stack-triangulation. Note that three different histories lead to the final triangulation. }
\end{figure}

The aim of this paper is to study these models of random maps. Particularly, we are interested in large maps when the number of faces tends to $+\infty$. It turns out that this model of triangulations is combinatorialy simpler that the set of all triangulations. Under the two probabilities $\qt_{2k}$ and $\ut_{2k}$ we exhibit a \sl global limit behavior \rm of these maps. \par

 A model of increasing quadrangulations is also treated at the end of the paper. In few words this model is as follows. Begin with the rooted square and successively choose a finite face $ABCD$, add inside a node $O$ and two new edges: $AO$ and $OC$ (or $BO$ and $OD$). When these two choices of pair of edges are allowed we get a model of quadrangulations that we were unable to treat as wanted (see Section \ref{qua-dur}). When only a suitable choice is possible, we get a model very similar to that of stack-triangulations that may be endowed also with two different natural probabilities.  The results obtained are, up to the normalizing constants, the same as those obtained for stack-triangulations.  For sake of briefness, only the case of stack-triangulations is treated in details.\medskip

We present below the content of the paper and a rough description of the results, the formal statements being given all along the paper.
\subsection{Contents}

In Section \ref{comb} we define formally the set of triangulations $\3{2n}$ and the probabilities $\ut_{2n}$ and $\qt_{2n}$. This section contains also a bijection between $\3{2n}$ and the set $\Tter_{3n-2}$ of ternary trees with $3n-2$ nodes deeply used in the paper.
In Section \ref{top} are presented the two topologies considered in the paper:
\begin{itemize}
\item[-] the first one is an ultra-metric topology called \sl topology of local convergence\rm. It aims to describe the limiting local behavior of a sequence of maps (or trees) around their roots,
\item[-] the second topology considered is the \sl Gromov-Hausdorff topology \rm on the set of compact metric spaces. It aims to describe the limiting behavior of maps (or trees) seen as metric spaces where the distance is the graph distance. The idea here is to normalize the distance in maps, say by their diameters, in order to observe a limiting behavior.
\end{itemize}
 In Section \ref{top} are also recalled some facts concerning Galton-Watson trees conditioned by the size, when the offspring distribution is $\nu_\ter=\frac13\delta_3+\frac23\delta_0$ (the tree is ternary in this case). In particular it is recalled that they converge under the topology of local convergence to an infinite branch, (the spine or infinite line of descent) on which are grafted some critical ternary Galton-Watson trees; 
rescaled by $n^{1/2}$  they converge for the Gromov-Hausdorff topology to the continuum random tree (CRT), introduced by Aldous \cite{ALD}.  \par
Section \ref{res} is devoted to the statements and the proofs of the main results of the paper concerning random triangulations  under $\ut_{2n}$, when $n\to+\infty$. The strongest theorems of this part, that may also be considered as the strongest results of the entire paper, are:
\begin{itemize}
\item[-] the weak convergence of  $\ut_{2n}$ for the topology of local convergence to a measure on infinite triangulations (Theorem \ref{loctri}),
\item[-] the weak convergence of the metric of stack-triangulations for the Gromov-Hausdorff topology (the distance being the graph distance divided by $\sqrt{6n}/11$) to the CRT (Theorem \ref{youp}). It is up to our knowledge, the only case where the convergence of the metric of a model of random maps is proved (apart from trees).
\end{itemize}

Section \ref{res2} is devoted to the study of $\3{2n}$ under $\qt_{2n}$. The behavior is very different from that under $\ut_{2n}$. First, there is no local convergence around the root, its degree going a.s. to $+\infty$. Theorem \ref{metconv} says that seen as metric spaces they converge normalized by $(6/11)\log n$, in the sense of the finite dimensional distributions, to the discrete distance on $[0,1]$ (the distance between different points is 1). Hence, there is no weak convergence for the Gromov-Hausdorff topology, the space $[0,1]$ under the discrete distance being not compact. Section \ref{azd} contains some elements stating the speed of growing of the maps (the evolution of the node degrees, or the size of a submap).\par
Section \ref{allquad} is devoted to the study of a model of quadrangulations very similar to that of stack-triangulations, and to some questions related to another family of growing quadrangulations.\par
Last, the Appendix, Section \ref{ap}, contains the proofs that have been extracted from the text for sake of clarity.

\subsection{Literature about stack-triangulations}
\rm
\label{lit}The fact that stack-triangulations are in bijection with ternary trees, used in this paper, seems to be classical and will be proved in Section \ref{yop}.\par 
Stack-triangulations appear in the literature for very various reasons. In Bernardi and Bonichon \cite{BB}, stack-triangulations are shown to be in bijection with intervals in the Kreweras lattice (and realizers being both minimal and maximal). The set of stack triangulations coincides also with the set of plane triangulations having a unique Schnyder wood (see Felsner and Zickfeld \cite{FZ}).\par
These triangulations appear also around the problem of graph uniquely 4-colorable.
A graph G is uniquely 4-colorable if it can be colored with 4 colors, and if every 4-coloring of G produces the same partition of the vertex set into 4 color classes. There is an old conjecture saying that the family of maps having this property is the set of stack-triangulations. We send the interested reader to B\"ohme \& al. \cite{BSV} and references therein for more information on the question.\par

\begin{figure}[htbp]\rm
\centerline{\includegraphics[height= 4.5 cm]{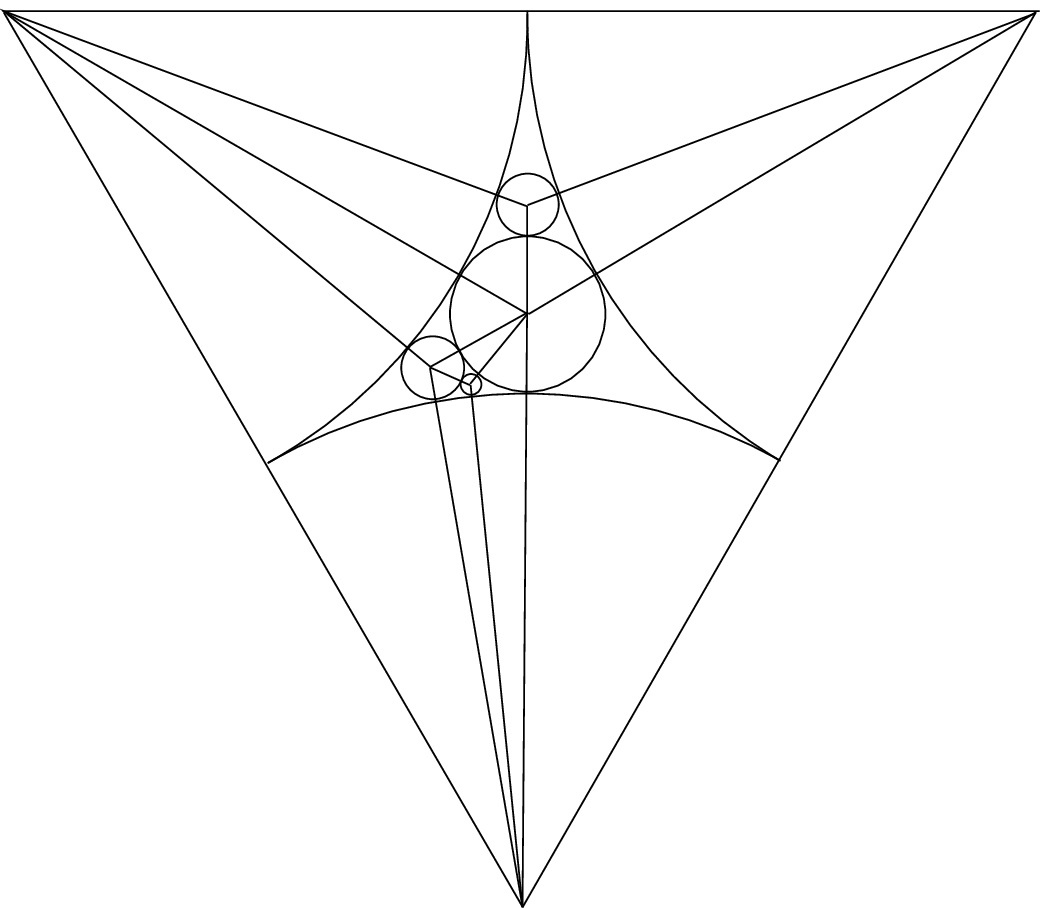}}
\captionn{\label{figapol}Construction of Apollonian's circles by successive insertions of circles (the starting point is three tangent circles). To get the triangulation associated, add an edge between the center of the new circle $C$ and the three centers of the circles tangent to $C$.}
\end{figure}

As illustrated on Figure \ref{figapol}, these triangulations appear also in relation with Apollonian circles. We refer to Graham \& al. \cite{GLMWY}, and to several other works of the same authors, for remarkable properties of these circles. \par
The so-called Apollonian networks, are obtained from Apollonian space-filling circles packing. First, we consider the Apollonian space-filling circles packing. Start with three adjacent circles as on Figure \ref{figapol}. The hole between them is filled by the unique circle that touches all three, forming then three new smaller holes. The associated triangulations is obtained by adding an edge between the center of the new circle $C$ and the three centers of the circles tangent to $C$.
If each time a unique hole receives a circle, the set of triangulation obtained are the stack-triangulations. If each hole received a circle all together, we get the model of Apollonian networks. We refer to Andrade \& al. \cite{AHAS} and references therein for some properties of this model of networks. \par 
The random Apollonian model of network studied by Zhou \& al. \cite{ZYW}, 
 Zhang \& al. \cite{ZRC},  and Zhang \& al. \cite{ZZ2} (when their parameters $d$ is 2) coincides with our model of stack-triangulations under $\qt$. Using physicist methodology and simulations they study among others the degree distribution (which is seen to respect a power-law) and the distance between two points taken at random (that is seen to be around $\log n$). \par

Darrasse and Soria \cite{DS} obtained the degree distribution on a model of ``Boltzmann'' stacked triangulations, where this time, the size of the quadrangulations is random, and uniformly distributed conditionally to its size. 

\medskip

We end the introduction by reviewing  the known asymptotic behaviors of quadrangulations and triangulations with $n$ faces under the uniform distribution (or close distributions in some sense). 

\subsection{Literature about convergence of maps} 

We refer to Angel \& Schramm \cite{AS}, Chassaing \& Schaeffer \cite{CS}  Bouttier \& al. \cite{BDG2} for an overview of the history of the study of maps from the combinatorial point of view, and to the references therein for the link with the 2-dimensional quantum gravity of physicists. We here focus on the main results concerning the convergence of maps. We exclude the results concerning trees (which are indeed also planar maps). \par
In the very last years, many studies concerning the behavior of large maps have been published. The aim in these works was mainly to define or to approach a notion of limiting map. Appeared then two different points of view, two different topologies to measure this convergence. \par
Angel \& Schramm \cite{AS} showed that the uniform distribution on the set of rooted triangulations with $n$ faces  (in fact several models of triangulations are investigated) converges weakly for a topology of local convergence (see Section \ref{tlconv}) to a distribution on the set of infinite but locally finite triangulations. In other words, for any $r$, the submap $S_r(n)$ obtained by keeping only the nodes and edges at distance smaller or equal to $r$ from the root vertex, converges in distribution toward a limiting random map $S_r$. By a theorem of Kolmogorov  this allows to show the convergence of the uniform measure on triangulations with $n$ faces to a measure on the set of infinite but locally finite rooted triangulations (see also Krikun \cite{MK} for a simple description of this measure). 
Chassaing \& Durhuus \cite{CD} obtained then a similar result, with a totally different approach, on uniform rooted quadrangulations with $n$ faces. \medskip

The second family of results concerns the convergence of rescaled maps: the first one in this direction has been obtained by Chassaing \& Schaeffer \cite{CS} who studied the limiting profile of quadrangulations.  The (cumulative) profile $(\Prof(k),k\geq 0)$ of a rooted graph, defined in Section \ref{asGH}, gives the successive number of nodes at distance smaller than $k$ from the root. Chassaing \& Schaeffer \cite[Corollary 4]{CS} showed that 
\[\l(\frac{\Prof((8n/9)^{1/4}x)}{n}\r)_{x\geq 0}\to \l({\cal J}[m,m+x]\r)_{x\geq 0}\]
where the convergence holds weakly in $D([0,+\infty),`R)$. The random probability measure ${\cal J}$ is  ISE the Integrated super Brownian excursion. ISE is the (random) occupation measure of the Brownian snake with lifetime process the normalized Brownian excursion, and $m$ is the minimum of the support of ${\cal J}$. The radius, i.e. the largest distance to the root, is also shown to  converge, divided by $(8n/9)^{1/4}$, to the range of ISE. Then,\\
-- Marckert \& Mokkadem \cite{MM2} showed the same result for pointed quadrangulations with $n$ faces,\\
-- Marckert  \& Miermont \cite{GM} showed that up to a normalizing constant, the same asymptotic holds for pointed rooted bipartite maps under Boltzmann distribution with $n$ faces, (the weight of a bipartite map is $\prod_{f \textrm{ face of m}}w_{\deg(f)}$ where the $(w_{2i})_{i\geq 0}$ is a ``critical sequence of weight''),  \\
-- Weill \cite{W} obtained the same results as those of \cite{GM} in the rooted case,\\
-- Miermont  \cite{GW2} provided the same asymptotics for rooted pointed Boltzmann maps  with $n$ faces with no restriction on the degree,\\
--  Weill and Miermont \cite{GW3}  obtained the same result as \cite{GW2} in the rooted case.\medskip

All these results imply that if one wants to find a (finite and non trivial) limiting object for rescaled maps, the edge-length in maps with $n$ faces has to be fixed to $n^{-1/4}$ instead of 1.  In Marckert \& Mokkadem \cite{MM2}, quadrangulations are shown to be obtained as the gluing of two trees, thanks to the Schaeffer's bijection (see \cite{S,CS,MM2}) between quadrangulations and well labeled trees. They introduce also a notion of random compact continuous map, ``the Brownian map'', a random metric space candidate to be the limit of rescaled quadrangulations. In \cite{MM2} the convergence of rescaled quadrangulations to the Brownian map is shown but not for a ``nice topology''. As a matter of fact, the convergence in \cite{MM2} is a convergence of the pair of trees that encodes the quadrangulations to a pair of random continuous trees, that also encodes, in a sense similar to the discrete case, a continuous object that they name the Brownian map. ``Unfortunately'' this convergence does not imply -- at least not in an evident way -- the convergence of the rescaled quadrangulations viewed as metric spaces to the Brownian map for the Gromov-Hausdorff topology. 

Some authors think that the Brownian map is indeed the limit, after rescaling, of classical families of maps (those studied in \cite{CS,MM2,GM,W,GW2,GW3}) for the Gromov-Hausdorff topology. Evidences in this direction have been obtained by Le Gall \cite{LGC} who proved the following result. He considers $M_n$ a $2p$-angulations with $n$ faces under the uniform law.  Then, he shows that at least along a suitable
subsequence, the metric space consisting of the set of vertices of $M_n$, equipped with the
graph distance rescaled by the factor $n^{1/4}$, converges in distribution as $n\to\infty$ towards
a limiting random compact metric space, in the sense of the Gromov-Hausdorff distance.
He proved that the topology of the limiting space is uniquely determined independently
of $p$ and of the subsequence, and that this space can be obtained as the quotient of
the CRT for an equivalence relation which is defined from Brownian
labels attached to the vertices. Then  Le Gall \& Paulin \cite{LGP} show that this limiting space is topologically a sphere. The description of the limiting space is a little bit different from the Brownian map but one may conjecture that these two spaces are identical. \par

Before coming back to our models and results we would like to stress on two points.\\
$\bullet$ The topology of local convergence (on non rescaled maps) and the Gromov-Hausdorff topology (on rescaled map) are somehow orthogonal topologies. The Gromov-Hausdorff topology considers only what is at the scaling size (the diameter, the distance between random points, but not the degree of the nodes for example). The topology of local convergence considers only what is at a finite distance from the root. In particular, it does not measure at all the phenomenons that are at the right scaling factor, if this scaling goes to $+\infty$. This entails that in principle one may not deduce any non-trivial limiting behavior for the Gromov-Hausdorff topology from the topology of local convergence, and vice versa. \\
$\bullet$ There is a conjecture saying that properly rescaled random planar maps conditioned to be large should converge to a limiting continuous random surface, whose law should not depend up to scaling constant from the family of reasonable maps that are sample. This conjecture still holds even if the family of stack-maps studied here converges to some objects that can not be the limit of uniform quadrangulations. The reason is that stack-maps are in some sense not reasonable maps.

\section{Stack-triangulations}
\label{comb}
\subsection{Planar maps}
A planar map $m$ is a proper embedding without edge crossing of a connected graph in the sphere. Two planar maps are identical if one of them can be mapped to the other by a homeomorphism that preserves the orientation of the sphere. A planar map is a quadrangulation if all its faces have degree four, and a triangulation if all its faces have degree three. There is a difference between the notions of planar maps and planar graphs, a planar graph having possibly  several non-homeomorphic embeddings on the sphere.   \par  
\begin{figure}[htbp]
\centerline{\includegraphics[height= 2.5 cm]{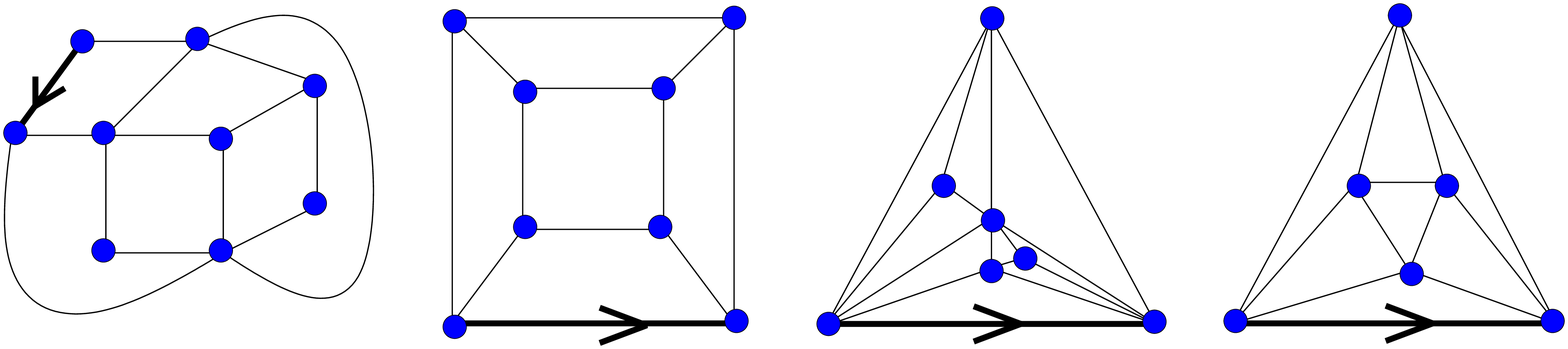}}
\captionn{\label{fig1}Two rooted quadrangulations and two rooted triangulations.}
\end{figure}
In this paper we deal with rooted planar maps $(m,E)$: an oriented edge $E=(E_0,E_1)$ of $m$ is distinguished. The point $E_0$ is called the root vertex of $m$. 
Two rooted maps are identical if the homeomorphism preserves also the distinguished oriented edge. Rooting maps like this allows to avoid non-trivial automorphisms. By a simple projection, rooted planar maps on the sphere are in one to one correspondence with rooted planar maps on the plane, where the root of the latter is adjacent to the infinite face (the unbounded face) and is oriented in such a way that the infinite face lies on its right, as on Figure \ref{fig1}. From now on, we work on the plane. \par

For any map $m$, we denote by $V(m),E(m), F(m), F^{\circ}(m)$ the sets of vertices, edges, faces and finite faces of $m$; for any $v$ in $V(m)$, we denote by $\deg(v)$ the degree of $v$. The graph distance $d_G$ between two vertices of a graph $G$ is the number of edges in a shortest path connecting them. The set of nodes of a map $m$ equipped with the graph distance denoted by $d_m$ is naturally a metric space. The study of the asymptotic behavior of $(m,d_m)$ under various distributions is the main aim of this paper.

\subsection{The stack-triangulations}
\label{def}

We build here $\3{2k}$ the set of \sl stack-triangulations \rm with $2k$ faces, for any $k\geq 1$. \par

\begin{figure}[htbp]\psfrag{x}{$x$}\psfrag{f}{$f$}
\centerline{\includegraphics[height= 2.2 cm]{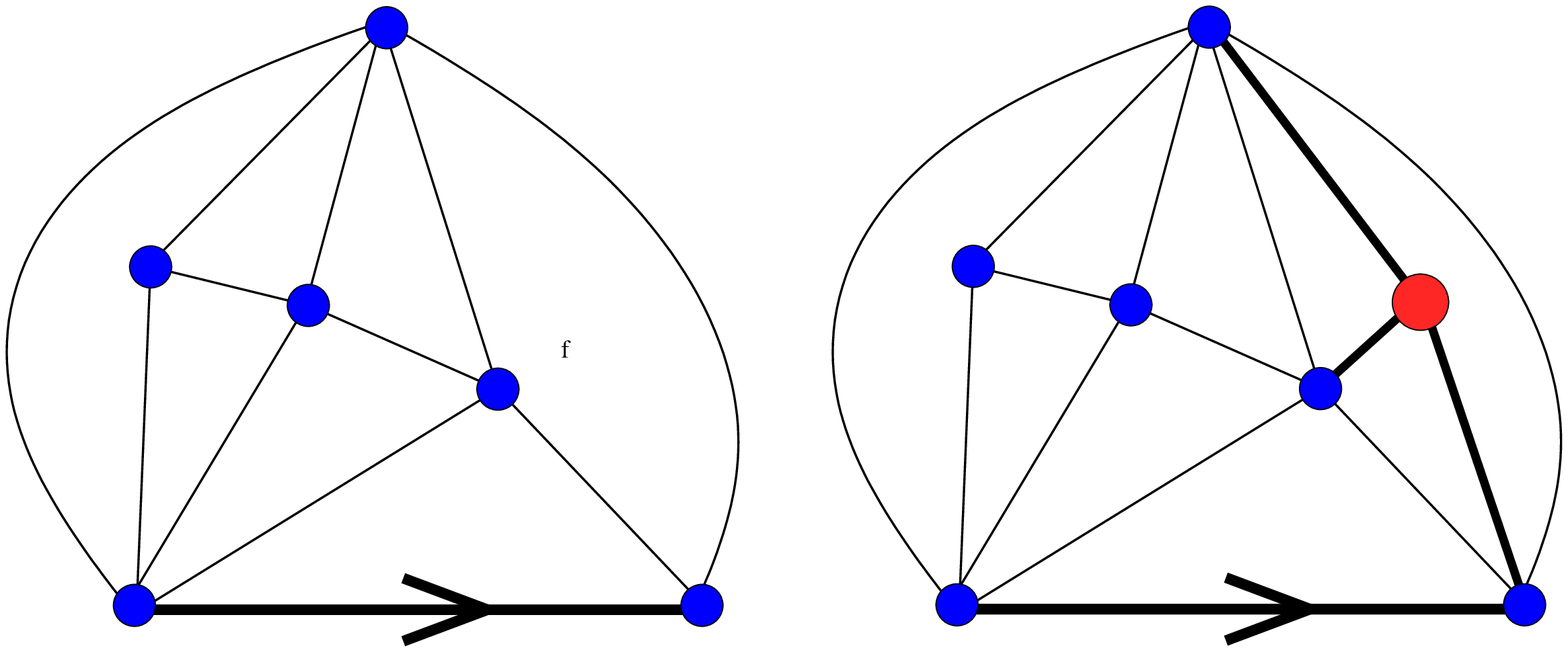}}
\captionn{\label{fig43}A triangulation $(m,f)$ with a distinguished face and its image by $\Phi$. }
\end{figure}

Set first $\3{2}=\{ \Theta\}$ where $\Theta$ denotes the unique rooted triangle (the first map in Figure \ref{fig2}). Assume that $\3{2k}$ is defined for some $k\geq 1$ and is a set of rooted triangulations with $2k$ faces. We now define $\3{2(k+1)}$. 
Let 
\[\3{2k}^{\bullet}=\{(m,f)~|~ m\in \3{2k},  f\in F^\circ(m)\}\] be the set of rooted triangulations from $\3{2k}$ with a distinguished finite face. 
We now introduce an application $\Phi$ from  $\3{2k}^{\bullet}$ onto the set of all rooted triangulations with $2(k+1)$ faces (we should write $\Phi_k$). 
For any $(m,f)\in \3{2k}^{\bullet}$, $\Phi(m,f)$ is the following rooted triangulation: draw $m$ in the plane, add a point $x$ inside the face $f$ and three non-crossing edges inside $f$ between $x$ and the three vertices of $f$ adjacent to $x$ (see Figure \ref{fig43}). The obtained map has $2k+2$ faces. \par
We call $\3{2(k+1)}=\Phi(\3{2k}^{\bullet})$ the image of this application.\medskip

On Figure \ref{fig1}, the first triangulation is in $\3{10}$ (see also Figure \ref{fig2}). The second one is not in $\3{8}$ since it has no internal node having degree 3.

\begin{defi}\label{history} We call internal vertex of a stack-triangulation $m$ every vertex of $m$ that is not adjacent to the infinite face (all the nodes but three).\\
We call history of a stack-triangulation $m_k\in\3{2k}$ any sequence $\big((m_i,f_i),i=1,\dots,k-1\big)$ such that $m_i\in\3{2i}$, $f_i\in F^\circ(m_i)$ and $m_{i+1}=\Phi(m_i,f_i)$. 
We let ${\cal H}(m)$ be the set of histories of $m$, and $H_\triangle(k)=\{{\cal H}(m) ~|~ m\in `3_{2k}\}$.
\end{defi}

We define here a special drawing ${\cal G}(m)$ of a stack-triangulation $m$. The embedding  ${\cal G}(\Theta)$ of the unique rooted triangle $\Theta$ is fixed at position $E_0=(0,0)$, $E_1=(1,0)$, $E_2=e^{i\pi/3}$ (where $E_0,E_1,E_2$ are the three vertices of $\Theta$, and $(E_0,E_1)$ its root). The drawing of its edges are straight lines drawn in the plane. To draw ${\cal G}(m)$ from ${\cal G}(m')$ when $m=\Phi(m',f')$, add a point $x$ in the center of mass of $f'$, and three straight lines between $x$ and the three vertices of $f'$ adjacent to $x$. The faces of ${\cal G}(m)$ hence obtained are geometrical triangles. Presented like this, ${\cal G}(m)$ seems to depend on the history of $m$ used in its construction, and thus we should have written ${\cal G}_h(m)$ instead of ${\cal G}(m)$, where the index $h$ would have stood for the history $h$ used. But it is easy to check (see Proposition \ref{yopp}) that if $h,h'$ are both in ${\cal H}(m)$ then ${\cal G}_{h'}(m)={\cal G}_{h}(m)$. 
\begin{defi}\label{cd}
The drawing ${\cal G}(m)$ is called the canonical drawing of $m$.
\end{defi}

\subsubsection{Two distributions on  $\3{2k}$}
\label{tw}
For any $k\geq 1$, we denote by $\ut_{2k}$ the uniform distribution on $\3{2k}$.\\
We now define a second probability $\qt_{2k}$. 
For this, we construct on a probability space $(\Omega,`P)$ a process $(M_n)_{n\geq 1}$ such that $M_n$ takes its values in $\3{2n}$ as follows:
first $M_1$ is the rooted triangle $\Theta$.
At time $k+1$, choose a finite face $F_{k}$ of $M_k$ uniformly among the finite faces of $M_k$ and this independently from the previous choices and set
\[M_{k+1}=\Phi(M_k,F_k).\]
We denote by $\qt_{2k}$ the distribution of $M_k$. Its support is exactly $\3{2k}$.

The weight of a map under $\qt_{2k}$ being proportional to its number of histories, 
it is easy to check that $\qt_{2k}\neq\ut_{2k}$ for $k\geq 4$.

\subsection{Combinatorial facts}
We begin this section where is presented the bijection between ternary trees and stack-triangulations with some considerations about trees.
\subsubsection{Definition of trees}
\label{deftree}
\begin{figure}[htbp]
\psfrag{e}{$\varnothing$}
\psfrag{11}{11}\psfrag{12}{12}\psfrag{211}{211}\psfrag{212}{212}
\psfrag{21}{21}\psfrag{13}{13}
\psfrag{2}{2}\psfrag{1}{1}
\centerline{\includegraphics[height= 2.5 cm]{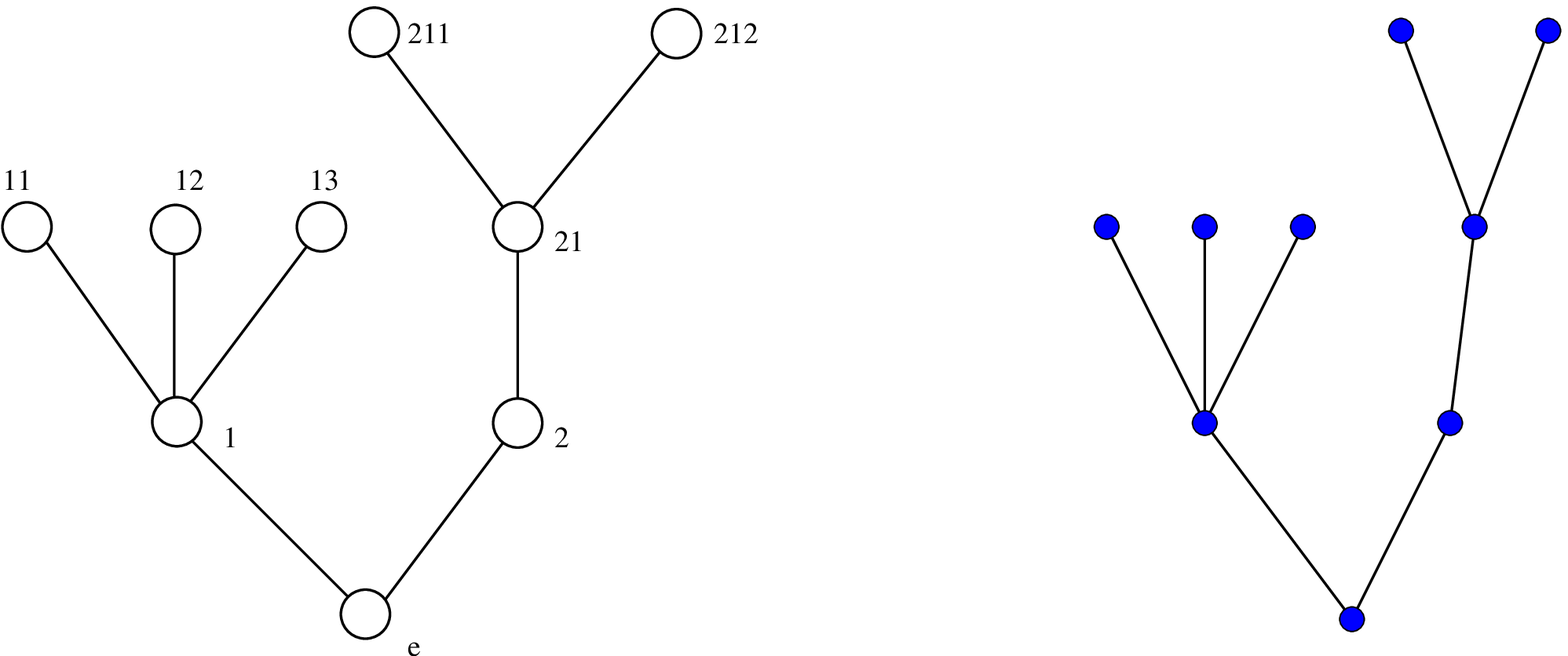}}
\captionn{\label{fig3}A rooted tree and its usual representation on the plane.   }
\end{figure}
Consider the set $W=\bigcup_{n\geq 0}\N^n$ of finite words on the alphabet $\mathbb{N}=\{1,2,3,\dots\}$ 
where by convention $\N^0=\{\varnothing\}$. For $u=u_1\ldots u_n,v=v_1\ldots v_m\in W$, we let 
$uv=u_1\ldots u_nv_1\ldots v_m$ be the concatenation of the words $u$ and $v$. 
\begin{defi}\label{tree}
A planar tree $\bt$ is a subset of $W$\\
$\bullet$ containing the root-vertex $\varnothing$,\\
$\bullet$ such that if $ui\in\bt$ for some $u\in W$ and $i\in\mathbb{N}$, then $u\in \bt$,\\
$\bullet$ and such that if  $ui\in\bt$ for some $u\in W$ and $i\in\mathbb{N}$, then $uj\in\bt$ for all $j\in\{1,\dots,i\}$. 
\end{defi}
We denote by $\T$ the set of planar trees.  For any $u\in \bt$, let  $c_u(\bt)=\max\{i~|~ui\in\bt\}$ be the number of \sl children \rm of $u$. The elements of a tree $\bt$ are called \sl nodes\rm, a node having no child a \sl leaf\rm, the other nodes the \sl internal nodes\rm.  The set of leaves of $t$ will be denoted by $\partial t$, and its set of internal nodes by $t^\circ$.\par A binary (resp. ternary) tree $\bt$ is a planar tree such that $c_u(\bt)\in\{0,2\}$ (resp. $c_u(\bt)\in\{0,3\}$) for any $u\in\bt$. We denote by $\Tbin$ and $\Tter$ the set of finite or infinite binary and ternary trees, and by $\Tbin_n$ and $\Tter_n$ the corresponding set of trees with $n$ nodes. \par

If $u$ and $v$ are two nodes in $\bt$, we denote by $u\wedge v$ the \sl deepest common ancestor \rm of $u$ and $v$, i.e. the largest word $w$ prefix to both $u$ and $v$ (the node $u\wedge v$ is in $\bt$). The length $|u|$ of a word $u\in W$ is called the height or depth  of $u$, or graph distance of $u$ to the root, if considered as a vertex of some tree. For $u=u_1\ldots u_n\in\bt$, we let 
$u[j]=u_1\ldots u_j$ and $[[\varnothing,u]]=\{\varnothing, u[1],\ldots, u[n]\}$ be the ancestral line of $u$ back to the root. For any tree $t$ and $u$ in $t$, the \sl fringe subtree \rm $t_u:=\{w~|~uw\in t\}$ is in some sense, the subtree of $t$ rooted in $u$. 
Finally recall that the lexicographical order (LO) on $W$ induces a total ordering of the nodes of any tree. 

\subsubsection{The fundamental bijection between stack-triangulations and ternary trees}
\label{fond-bij}

Before explaining the bijection we use between $\3{2K}$ and $\Tter_{3K-2}$ we define a function $\Gamma$ which will play an eminent role in our asymptotic results concerning the metrics in maps. 
Let $W_{1,2,3}$ be the set of words containing at least one occurrence of each 
element of $\Sigma_3=\{1,2,3\}$ as for example  $321$, $123$, $113211213123$.
Let $u=u_1\dots u_k$ be a word on the alphabet $\Sigma_3$. Define 
 $\tau_1(u):=0$ and $\tau_2(u) := \inf \{i, i>0 , u_i=1\}$, the rank of
the first apparition of $1$ in $u$. For $j\ge3$,  define
\[\tau_j(u) :=\inf\{i~|~ i > \tau_{j-1}(u) \text{ such that }
u_{1+\tau_{j-1}(u)}\dots u_i\in W_{1,2,3} \}.\]
This amounts to decomposing $u$ into subwords, the first one ending when the first 1 appears, the subsequent ones ending each time that every of the three letters 1, 2 and 3 have appeared again. For example if $u=22123122131$ then $\tau_{1}(u)=0,\tau_2(u)=3,\tau_3(u)=6,\tau_4(u)=10$. Denote by 

\begin{equation}
\Gamma(u)=\max\{i~|~ \tau_i(u)\leq |u|\}
\end{equation} the number of these non-overlapping 
subwords.
Further for two words (or nodes) $u=w a_1 \dots a_k$ and $v=w b_1 \dots b_l$ with $a_1\neq b_1$ and $w=u\wedge v$, set 
\begin{equation}
\Gamma(u,v)= \Gamma(a_1\dots a_k)+\Gamma(b_1\dots b_l).
\end{equation}
We call the one or two parameters function $\Gamma$ the \sl passage function\rm.
We know describe a bijection $\Psi_K^{`3}$ between $\3{2K}$ and $\Tter_{3K-2}$ having a lot of important properties.
\begin{pro}\label{yop}
For any $K\geq 1$ there exists a bijection 
\[\app{\Psi_K^{`3}}{\3{2K}}{\Tter_{3K-2}}{m}{t:=\Psi_K^{`3}(m)}\]
such that:\\
$(i)$ $(a)$ Each internal node $u$ of $m$ corresponds bijectively to an internal node $v$ of $t$. We denote for sake of simplicity by $u'$ the image of $u$.\\
$(b)$ Each leaf of $t$ corresponds bijectively to a finite triangular face of $m$. \\
$(ii)$ For any $u$ internal node of $m$, $\Gamma(u')=d_m(root,u).$ \\
$(ii')$ For any $u$ and $v$ internal nodes of $m$
\begin{equation}\label{ezqus}
\l|d_{m}(u,v)-\Gamma(u',v')\r|\leq 4.
\end{equation}
\noindent$(iii)$ Let $u$ be an internal node of $m$. We have
\[\deg_m(u)=\#\{v'\in t^\circ~|~ v'=u'w', w'\in 1L^\star_{2,3}\cup 3L^\star_{1,2}\cup 2L^\star_{1,3}\},\]
where $\{v'\in t^\circ~|~ v'=u'w', w'\in 1L^\star_{2,3}\cup 2L^\star_{1,3}\cup 3L^\star_{1,2}\}$ is the union of the subtrees of $t^\circ$ rooted in $u'1$, $u'2$ and $u'3$ formed by the ``binary trees'' having no nodes containing a 1, resp. a 2, resp a 3.
\end{pro}
We will write $\Psi^\triangle$ instead of $\Psi^\triangle_K$ when no confusion on $K$ is possible.\par
The last property of the Proposition \ref{yop} can be found in Darasse \& Soria \cite{DS}; we give below a proof for the reader convenience. The quote around binary trees signal that by construction these binary tree like structures do not satisfy the requirement of Definition \ref{tree}.

\begin{note}\rm\label{aut-bij}
The existence of a bijection between  $\3{2K}$ and $\Tter_{3K-2}$ follows the ternary decomposition of the maps in $\3{2K}$, as illustrated on Figure \ref{fig:dec_tri}: in the first step of the construction of $m$, the insertion of the three first edges incident to the node $x$ in the triangle $\Theta$ splits it into three parts that behave clearly as stack-triangulations (see Section \ref{def}). The node $x$ may be recovered at any time since it is the unique vertex incident to the three vertices incident to the infinite face. 
\begin{figure}[htbp]
\psfrag{A}{$E_0$}
\psfrag{a}{$\3{2k}(m)=$}
\psfrag{b}{$\ttt(m_1)$}
\psfrag{c}{$\ttt(m_2)$}
\psfrag{d}{$\ttt(m_3)$}
\psfrag{B}{$E_1$}
\psfrag{C}{$E_2$}
\psfrag{O}{$O$}\psfrag{x}{$x$}
\psfrag{O'}{$O'=\varnothing$}
\psfrag{O2}{$O_2$}
\psfrag{O3}{$O_3$}
\psfrag{m}{$m$}
\psfrag{m1}{$m_1$}
\psfrag{m2}{$m_2$}
\psfrag{m3}{$m_3$}
\centerline{\includegraphics[height= 4.5 cm]{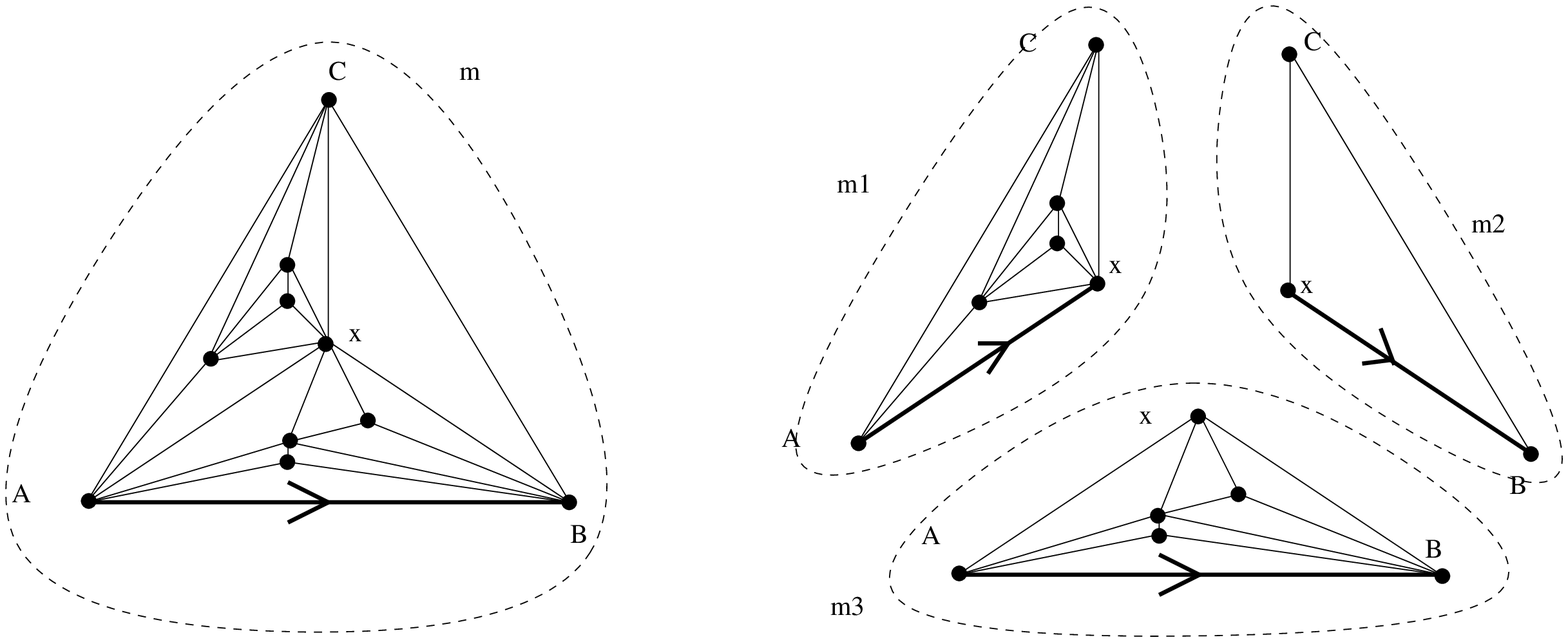}}
\captionn{\label{fig:dec_tri} Decomposition of a stack-triangulation using the recovering of the first inserted node.}
\end{figure}
\end{note}

Proposition \ref{yop} follows readily Proposition \ref{yopp} below, which is a history-dependent analogous. 
We denote by $\Tterb_{3n+1}:=\{(t,u)~|~ t\in \Tter_{3n+1}, u\in \partial t\}$ the set of ternary trees with $3n+1$ nodes with a distinguished leaf. Very similarly with the function $\Phi$ defined in Section \ref{def}, we define the application $\phi$ from $\Tterb_{3k+1}$ into $\Tter_{3k+4}$ as follows; for any $(t,u)\in \Tterb_{3k+1}$, let $t':=\phi(t,u)$ be the tree $t\cup\{u1,u2,u3\}$ obtained from $t$ by the replacement of the leaf $u$ by an internal node having 3 children. As for maps (see Definition \ref{history}), for any tree $t\in \Tter_{3k-2}$, a history of $t$ is a sequence $h'=\big((t_i,u_i),i=1,\dots,k-1\big)$ such that $(t_i,u_i)\in \Tterb_{3i-2}$ and  $t_{i+1}=\phi(t_i,u_i)$. The set of histories of $t$ is denoted by ${\cal H}(t)$, and we denote $H_{\cal T}(k)=\{{\cal H}(t)~|~ t\in \Tter_{3k-2}\}$. 

Notice that if $h=\big((m_i,f_i),i=1,\dots,K-1\big)$ is a history of $m$ then for any $j\leq K$, $h_j:=\big((m_i,f_i),i=1,\dots,j-1\big)$ is a history of $m_j$ (and the same property holds true for the histories $h'_j:=\big((t_i,u_i),i=1,\dots,j-1\big)$ associated with $h'=\big((t_i,u_i),i=1,\dots,k-1\big)$).

\begin{pro}
\label{yopp}
For any $K\geq 1$ there exists a bijection 
\[\app{\psi_K^{`3}}{H_\triangle(K)}{H_{\cal T}(K)}{h}{h'}\]
such that 
:\\
$(i)$ The family $(\psi_K^{`3},K\geq 1)$ is consistent: if $\psi_K^{`3}(h)=h'$ then for any $j\leq K$,
\[\psi_j^{`3}(h_j)=h'_j;\]
in other words $m_{i+1}=\Phi(m_i,f_i)$ is sent on $t_{i+1}:=\phi(t_i,u_i).$\\
$(i')$ Robustness: $h_1$ and $h_2$ are two histories of $m$ iff $\psi_K^{`3}(h_1)$ and $\psi_K^{`3}(h_2)$ are histories of the same tree.\\
$(ii)$ Each finite face of $m$ corresponds to a leaf $u$ of $t$,
and each internal node $u$ of $m$ corresponds bijectively to an internal node $v$ of $t$. We denote for sake of simplicity by $u'$ the image of $u$.\\
$(ii)$ For any $u$ internal node of $m$, $\Gamma(u')=d_m(root,u).$ \\
$(ii')$ For any $u$ and $v$ internal nodes of $m$
\begin{equation}\label{ezquss}
\l|d_{m}(u,v)-\Gamma(u',v')\r|\leq 4.
\end{equation}
\noindent$(iii)$ Let $u$ be an internal node of $m$. We have
\[\deg_m(u)=\#\{v'\in t^\circ~|~ v'=u'w', w'\in 1L^\star_{2,3}\cup 3L^\star_{1,2}\cup 2L^\star_{1,3}\},\]
where $\{v'\in t^{\circ}~|~ v'=u'w', w'\in 1L^\star_{2,3}\cup 2L^\star_{1,3}\cup 3L^\star_{1,2}\}$ is the union of the subtrees of $t$ rooted in $u'1$, $u'2$ and $u'3$ formed by the ``binary trees'' having no nodes containing a 1, resp. a 2, resp a 3. 
\end{pro}
\proof The proof of $(i), (i'), (ii), (ii')$ follows directly the explicit construction presented below and the property that during the construction of stack triangulation, the insertion in a given face $f$ does not modify the other faces (as well as in a tree, the insertions in a subtree does not change the other subtrees). In one word, this construction raises on a canonical association of a triangular face with a node of $W$. We stress on this point, and we recall the canonical drawing of the Definition  \ref{cd}: thanks to the canonical drawing there is a sense to talk of a face $f$ without referring to a map, and thanks to our construction of trees, there is a sense to talk of a node $u$ (which is a word) without referring to a tree. We will call \sl canonical face \rm a geometrical face corresponding to a canonical drawing.
The bijection $\psi$ presented below works at this level and associates with a canonical face a word and thus has immediately the properties of consistence and robustness stated in the Proposition.  We proceed as follows.\par

In the case  $K=1$, $\3{2K}=\3{2}$ contains only the rooted triangle $\Theta$, and $\Tter_{3K-2}=\Tter_1$ contains only the tree  reduced to the root $\varnothing$, which is a leaf in this case: hence set $\psi^{`3}_2(\Theta)=\{\varnothing\}$.\par

We have mainly to explain how ``canonically'' associate with a growing of a map a suitable growing of the corresponding trees. Assume then that each finite canonical face (of each maps of $\triangle_{2K}$) are associated with some nodes (of the trees of $\Tter_{3K-2}$): that is, if the canonical drawing of the maps is used, whatever the maps is considered, a given face is sent on the same word. We then  write $\psi^\triangle(f)=u$ for a face $f$ and a word $u$.
We associate canonically with each finite face $f$ of any map $m$ an oriented edge $O(f)$ as follows. At first $\Theta$, the rooted triangle has three vertices $E_0,E_1,E_2$, $(E_0,E_1)$ being the root. For this finite face, choose $O(f)=(E_0,E_1)$. 
Assume now that each finite face $f$ of $m(i)$ owns an oriented edge $O(f)$.
Assume that the face $f(i)=(A,B,C)$ has for associated oriented edge
$O(f(i))=(A,B)$, and that it is sent on the node $u$. We set as respective associated edges for the three new faces $(B,C,x)$, $(A,x,C)$, $(A,B,x)$ the edges  $(A,x)$, $(x,B)$, and $(A,B)$ (they are oriented in such a way that the infinite face lies on the right of the faces seen as maps, and allow a successive decomposition). We associate respectively to these faces the nodes  $u1$ , $u2$ and  $u3$. This construction is indeed canonical, in the sense that if a face belongs to two canonical drawings ${\cal G}(m)$ and ${\cal G}(m')$, then the new faces obtained after growing are sent by our construction to the same nodes.

\begin{figure}[htbp]
\psfrag{a}{$A$}
\psfrag{b}{$B$}
\psfrag{c}{$C$}
\psfrag{x}{$x$}
\centerline{\includegraphics[height= 2.5 cm]{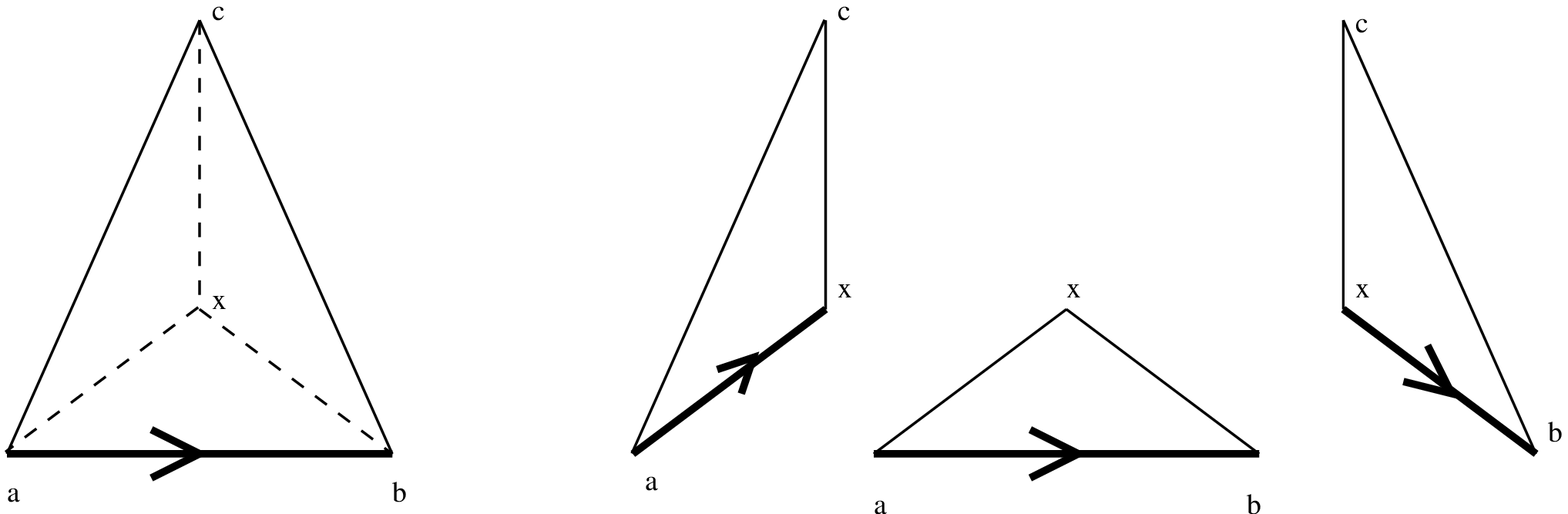}}
\captionn{\label{st} Heritage of the canonical orientation of the faces. If the first face is sent on $u$, then the other ones, from left to right are sent on $u2, u3$ and $u1$} 
\end{figure}
This leads easily, by recurrence on $K$, to $(i), (i'), (ii)$ and $(ii')$. \par

In order to prove the other assertions of the Proposition, we introduce the notion of type of a face, and of a node. For any face $(u,v,w)$ in $m$, define 
\begin{equation}
\type(u,v,w):=\l(d_m(E_0,u),d_m(E_0,v),d_m(E_0,w)\r),
\end{equation}
the distance of $u,v,w$ to the root-vertex of $m$. 
Since $u$, $v$, and $w$ are neighbors, the type of any triangle is $(i,i,i)$, $(i,i,i+1)$, $(i,i+1,i+1)$ for some $i$, or a permutation of this. The types of the faces arising in the construction of $m$ are also well defined, since the insertions do not change the distance between the existing nodes and the root. 
We then prolong the construction of $\Phi$ given above, and mark the nodes of $t$ with the types of the corresponding faces. For any internal node $u'\in t$ with $\type(u)=(i,j,k)$,
\begin{equation}\label{evolo}
\left\{
\begin{array}{ccccl}
\type(u1)=(&1+i\wedge j \wedge k,&j,&k&),\\
\type(u2)=(&i,&1+i\wedge j \wedge k,&k&),\\
\type(u3)=(&i,&j,&1+i\wedge j \wedge k&)
\end{array}\right.
\end{equation}
as one can easily check with a simple figure: this corresponds as said above to the fact that if the leaf $u$ is associated with the ``empty'' triangle $(A,B,C)$, then the insertion of a node $x$ in $(A,B,C)$ is translated by the insertion in the tree of the nodes $u1$ (resp. $u2$, $u3$) associated with $(x,B,C)$ (resp. $(A,x,B)$, $(A,B,x)$). Formula \eref{evolo} gives then the types of these three faces.
Using that $\type(\varnothing)=(0,1,1),$ giving $t$ the types of all nodes are known and are obtained via the deterministic evolution rules \eref{evolo}. 
\begin{figure}[htbp]
\psfrag{A}{$E_0$}
\psfrag{B}{$E_1$}
\psfrag{C}{$E_2$}
\centerline{\includegraphics[height= 4 cm]{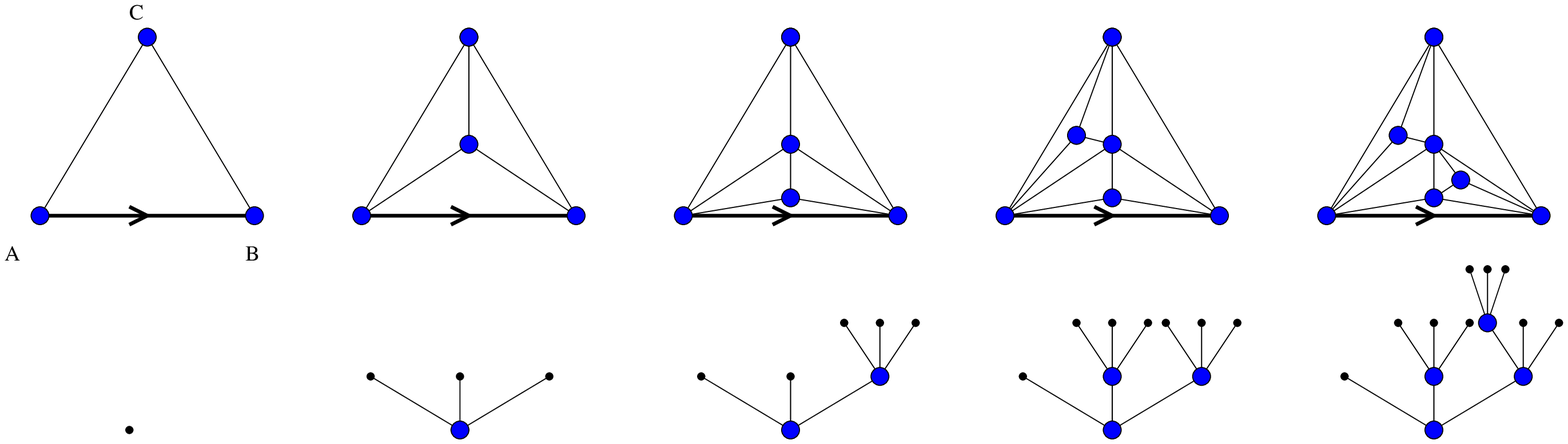}}
\captionn{\label{fig:tree_tri} Construction of the ternary tree associated with an
history of a stack-triangulation} 
\end{figure}
The distance of any internal node $u$ to the root of $m$ is computed as follows:  assume that 
$u$ has been inserted at a certain date in a face $f=(A,B,C)$.
Then clearly its distance to the root vertex is
\[d_m(E_0,u)=g(type (f)),\]
where $g(i,j,k)= 1+(i\wedge j \wedge k).$
Moreover, since an internal node in $m$ corresponds to the insertion of three children in the tree, each internal node $u$ of $m$ corresponds to an internal node $u'$ of $t$ 
and  
\[d_m(E_0,u)=g(\type(u')).\]

It remains to check that for any $u'\in t$, 
\begin{equation}
g(\type(u'))= \Gamma(u')
\end{equation} as defined above. This is a simple exercise: the initial type (that of $\varnothing$) varies along a branch of $t$ only when a 1 occurs in the nodes. Then the type passes from $(i,i,i)$ to $(i+1,i+1,i+1)$ when the three letters 1, 2 and 3 have appeared: this corresponds to the incrementation of the distance to the root in the triangulation.\\
$(ii')$ Consider $u$ and $v$ two internal nodes of $m$. The node $w'=u'\wedge v'$ corresponds to the smallest triangle $f$ containing $u$ and $v$. For some $a\neq b,$ belonging to $\{1,2,3\}$, the nodes $wa$ and $wb$ correspond to two triangles  containing respectively $u$ and $v$ (but not both). It follows that the distance $d_m(u,v)$ is equal up to 2, to $d_m(u,w)+d_m(w,v)$. Let us investigate now the relation between $w$ and $u$ and $\Gamma(a_1\dots a_j)$ in the case where $u=wa_1\dots a_j$. Each triangle appearing in the construction of $m$ behaves as a copy of $m$ except that its type is not necessarily $(i,i+1,i+1)$ (as was the type of $\varnothing$). Then the distance of the node $u=wa_1\dots a_j$ to $w$ may be not exactly $\Gamma(a_1\dots a_j)$. We now show that
\[\l|d_m(w,u)-\Gamma(a_1\dots a_j)\r|\leq 1.\]
This difference comes from the initialization of the counting of the non-overlapping subwords from $W_{1,2,3}$ in $a_1\dots a_j$. In the definition of $\Gamma$,  $\tau_2$ has a description different from the other $\tau$ since only a 1 is needed to reach a face having type $(i,i,i)$, the face $(1,1,1)$. Here again, according to the type of the face of $u$, the first $\tau$ corresponding to the reaching time of a face of type $(i,i,i)$ may have a different form than $\tau_2$ defined at the beginning of the Section: this $\tau$ is the waiting time of a 1, a 2, a 3 or of two letters among $\{1,2,3\}$. In any case, the corresponding $\Gamma$, say $\Gamma'$ verifies clearly
\[|\Gamma(z)-\Gamma'(z)|\leq 1\]
for any word $z$ on the alphabet $\Sigma_3$. \par
$(iii)$ This follows the description given in $(ii)$, since the degree of a node $u$ is the number of nodes at distance 1 of $u$. ~$\hfill\Box$

\subsection{Induced distribution on the set of ternary trees}
\label{ind-dist}
The bijection $\Psi_K^{`3}$ transports obviously the distributions $\ut_{2K}$ and $\qt_{2K}$ on the set of ternary trees $\Tter_{3K-2}$. \\
1) First, the distribution 
\begin{equation}\label{uniter}
\uter_{3K-2}:=\ut_{2K}\circ (\Psi_K^{`3})^{-1}
\end{equation}
is simply the uniform distribution on $\Tter_{3K-2}$ since $\Psi_K^{`3}$ is a bijection.\\
2) The distribution
\begin{equation}\label{hister}
\qter_{3K-2}:=\qt_{2K}\circ (\Psi_K^{`3})^{-1}
\end{equation}
is the distribution giving a weight to a tree proportional to the number of histories of the corresponding triangulation.\par

We want to give here another representation of the distribution $\qter_{3K-2}$.

\begin{defi} We call increasing ternary tree $\bt=(T,l)$ a pair such that:\\
$\bullet$ $T$ is the set of internal nodes of a ternary tree, \\
$\bullet$ $l$ is a bijective application between $T$ (viewed as a set of nodes) onto $\{1,\dots,|T|\}$ such that $l$ is increasing along the branches (thus $l(\varnothing)=1$).
\end{defi}
Notice that $T$ is not necessarily a tree as defined in Section \ref{deftree}: for example $T$ may be $\{\varnothing, 2\}$.

Let ${\cal I}_K^\ter$ denotes the set of increasing ternary trees $(T,l)$ such that  $|T|=K$ (i.e. $T$ is the set of internal nodes of a tree in  $\Tter_{3K+1}$).\par
The number of histories of a ternary tree $t\in \Tter_{3K-2}$ is  given by the 
\[w_{K-1}(t^\circ)=\#\{(t^\circ,l)\in {\cal I}_{K-1}^\ter\}\]
the number of increasing trees having $t^\circ$ as first coordinate, in other words, with shape $t^\circ$. Indeed, in order to record the number of histories of $t$ an idea is mark the internal nodes of $t$ by their apparition time (the root is then marked 1). Hence the marks are increasing along the branches, and there is a bijection between $\{1,\dots,K-1\}$ and the set of internal nodes of $t$. Conversely, any labeling of $t^\circ$ with marks having these properties corresponds indeed to a history of $m$. Thus
\begin{lem}For any $K\geq 1$, the distribution $\qter_{3K-2}$ has the following representation: for any $t\in\Tter_{3K-2}$,
\[\qter_{3K-2}(t)=C_{K-1}\cdot w_{K-1}(t^\circ)\]
where $C_{K-1}$ is the constant $C_{K-1}:=\l(\sum_{t'\in \Tter_{3K-2}} w_{K-1}({t'}^\circ)\r)^{-1}.$
\end{lem}

\section{Topologies}
\label{top}
\subsection{Topology of local convergence}
\label{tlconv}
The topology induced by the distance $d_L$ defined below will be called ``topology of local convergence''. 
Its aim is to describe an asymptotic behavior of maps (or more generally graphs) around their root. We stress on the fact that the limiting behavior is given under no rescaling.\par
We borrow some considerations from Angel \& Schramm \cite{AS}.  Let ${\cal M}$ be the set of rooted maps $(m,e)$ where $e=(e_0,e_1)$ is the distinguished edge of $m$. The maps from ${\cal M}$ are not assumed to be finite, but only locally finite, i.e. the degree of the vertices are finite.
For any $r\geq 0$, denote by $B_m(r)$ the map having as set of vertices 
\[V(B_m(r))=\{u\in V(m)~|~ d_m(u,e_0)\leq r\},\]
the vertices in $m$ with graph distance to $e_0$ non greater than $r$, and having as set of edges, the edges in $E(m)$ between the vertices of $V(B_m(r))$. 

For any $m=({\sf m}_1,e)$ and $m'=({\sf m}',e')$ in ${\cal M}$ set
\begin{equation}\label{dl}
d_L(m,m')={1}/({1+k})
\end{equation}
where $k$ is the supremum of the radius $r$ such that $B_m(r)$ and $B_{m'}(r)$ are equals as rooted maps. The application $d_L$ is a metric on the space ${\cal M}$.
A sequence of rooted maps converges to a given rooted map $m$ (for the metric $d_L$) if eventually they are equivalent with $m$ on arbitrarily large combinatorial balls around their root. In this topology, all finite maps are isolated points, and infinite maps are their accumulation points.
The space ${\cal M}$ is complete for the distance $d_L$ since given a Cauchy sequence of locally finite embedded rooted maps it is easy to see that it is possible to choose for them embeddings that eventually agree on the balls
of any fixed radius around the root. Thus, the limit of the sequence exists
(as a locally finite embedded maps). In other words, the space $\T$ of (locally finite embedded rooted) maps is complete.

The  space of triangulations (or of quadrangulations) endowed with this metric is not compact since it is easy to find a sequence of triangulations being pairwise at distance 1. The topology on the space of triangulations induces a weak topology on the linear space of measures supported on planar triangulations. 

\subsection{Gromov-Hausdorff topology}
\label{GHT}
The other topology we are interested in will be the suitable tool to describe the convergence of rescaled maps to a limiting object. The point of view here, is to consider maps endowed with the graph distance as metric spaces. The topology considered -- called the Gromov-Hausdorff topology -- is the topology of the convergence of compact (rooted) metric spaces. We borrow some considerations from Le Gall \& Paulin \cite{LGP} and from Le Gall \cite[Section 2]{LGC2}. We send the interested reader to these works and references therein.

First, recall that the Hausdorff distance in a metric space $(E,d_E)$ is a distance between the compact sets of $E$; for $K_1$ and $K_2$ compacts in $E$, 
\[d_{Haus(E)}(K_1,K_2)=\inf\{r ~|~ K_1 \subset K_2^r, K_2 \subset K_1^r\}\]
where $K^r=\cup_{x\in K} B_E(x,r)$ is the union of open balls of radius $r$ centered on the points of $K$. Now, given two pointed(i.e. with a distinguished node) compact metric spaces $((E_1,v_1),d_1)$ and $((E_2,v_2),d_2)$, the Gromov-Hausdorff distance between them is
\[d_{GH}(E_1,E_2)=\inf d_{Haus(E)}(\phi_1(E_1),\phi_2(E_2))\vee d_E(\phi_1(v_1),\phi_2(v_2))\,\]
where the infimum is taken on all metric spaces $E$ and all isometric embeddings $\phi_1$ and $\phi_2$ from $(E_1,d_1)$ and $(E_2,d_2)$ in $(E,d_E)$. Let $\mathbb{K}$ be the set of all isometric classes of compact metric spaces, endowed with the Gromov-Hausdorff distance $d_{GH}$. It turns out that $(\mathbb{K},d_{GH})$ is a complete metric space, which makes it appropriate to study the convergence in distribution of $\mathbb{K}$-valued random variables.
Hence, if $(E_n,d_n)$ is a sequence of metric spaces, $(E_n,d_n)$ converges for the Gromov-Hausdorff topology if there exists a metric space $(E_{\infty},d_{\infty})$ such that $d_{GH}(E_n,E_{\infty})\to 0$. \par
The Gromov-Hausdorff convergence is then a consequence of any convergence of $E'_n$ to $E'_{\infty}$, when $E'_n$ and $E'_{\infty}$ are some isomorphic embeddings of $E_n$ and $E_{\infty}$ in a common metric space $(E,d_E)$. In the proofs, we exhibit a space $(E,d_E)$ where this convergence holds; hence, the results of convergence we get are stronger than the only convergence for the Gromov-Hausdorff topology. In fact, it holds for a sequence of parametrized spaces.

\subsection{Galton-Watson trees conditioned by the size}
\label{reo}
Consider  $\nu_\ter:=\frac23\delta_0+\frac13\delta_3$ as a (critical) offspring distribution of a Galton-Watson (GW) process starting from one individual.
Denote by $P^{\ter}$ the law of the corresponding GW family tree; we will also write $P^{\ter}_n$ instead of $P^{\ter}\l(~.~ \big |\, |\bt|=n\r)$.
\begin{lem}\label{c1} $P^{\ter}_{3n+1}$ is the uniform distribution on $\Tter_{3n+1}$.
\end{lem}
\proof A ternary tree $t$ with $3n+1$ nodes has $n$ internal nodes having degree 3 and $2n+1$ leaves with degree 0. Hence $P^{\ter}_{3n+1}(\{t\} )=3^{-n}(2/3)^{2n+1}/P^{\ter}(\Tter_{3n+1})$.
This is constant on $\Tter_{3n+1}$ and has support $\Tter_{3n+1}$.~$\hfill\Box$

The conclusion is that for any $K\geq 1$
\begin{equation}\label{equality}
P^{\ter}_{3K-2}=\uter_{3K-2}.
\end{equation}
Following \eref{uniter}, this gives us a representation of the uniform distribution on $\Tter_{3K-2}$ in terms of conditioned GW trees. This will be our point of view in the sequel of the paper.\medskip

The asymptotic behavior of GW trees under $`P_n$ is very well studied. We focus in this section on the limiting behavior 
under the Gromov-Hausdorff topology and the topology of local convergence. The facts described here will be used later in the proof of the theorems stating the convergence of stack-triangulations. In addition we stress on the fact that the limit of rescaled stack-maps under the uniform distribution is the same limit as the one of GW trees: the continuum random tree. The rest of this Section is taken from the existing literature and is given for the reader convenience.

\subsubsection{Local convergence of Galton-Watson trees conditioned by the size}
\label{spine}
We endow ${\Tter}$ with the local distance $d_L$ defined in \eref{dl}. Under this metric, the accumulation points of sequences of trees $(t_K)$ such that $|t_K|=3K-2$ are infinite trees. 
It is known that the sequence $(P^{\ter}_{3K-2})$ converges weakly for the topology of local convergence. Let us describe a random tree $\bt_{\infty}^\ter$ under the limit distribution, denoted by $P^{\ter}_{\infty}$.\par
 Let ${\cal W}_3$ be the infinite ternary tree containing all words on  $\Sigma_3=\{1,2,3\}$ 
and let $(X_i)$ be a sequence of i.i.d. r.v. uniformly distributed on $\Sigma_3$.
Define
\begin{equation}
L^\ter_{\infty}=(X(j), j\geq 0)
\end{equation}
the infinite path in ${\cal W}_3$ starting from the root ($\varnothing$) and containing the words $X(j):=X_1\dots X_j$ for any $j\geq 1$. Take a sequence $(t(i))$ of GW trees under $P^\ter$ and graft them 
on the neighbors of $L_{\infty}^\ter$, that is on the nodes of ${\cal W}_3$ at distance 1 of $L_{\infty}^\ter$ (sorted according to the LO). The tree obtained is  $\bt_{\infty}^\ter$. 
In the literature the branch $L_\infty^\ter$ is called the \sl spine \rm or the \sl infinite line of descent \rm in $\bt_{\infty}^\ter$.\par
\begin{pro}(Gillet \cite{FG})\label{loctree}
When $n\to+\infty$, $P^{\ter}_{3n+1}$ converges weakly to $P^{\ter}_{\infty}$ for the topology of local convergence. 
\end{pro}
This result is due to Gillet \cite[Section III]{FG} (see Theorems III.3.1, III.4.2, III.4.3, III.4.4). 
\begin{note}\rm The distribution $P^{\ter}_{\infty}$ is usually called ``size biased GW trees''. We send the interested reader to Section 2 in Lyons \& al. \cite{LPP} to have an overview of this object. In particular, this distribution is known to be the limit of critical GW trees conditioned by the non extinction.
\end{note}

\subsubsection{Gromov-Hausdorff convergence of rescaled GW trees}
\label{ghg}
We present here the limit of rescaled GW trees conditioned by the size for the Gromov-Hausdorff topology. We borrow some considerations from Le Gall \& Weill \cite{legweill} and Le Gall \cite{LGC2}. \par
\begin{figure}[htbp]
\psfrag{0}{0}\psfrag{x}{$x$}
\psfrag{1}{1}
\psfrag{s}{$s$}
\psfrag{s'}{$s'$}
\psfrag{t}{$t$}
\psfrag{f(s)}{$f(s)$}
\psfrag{f(t)}{$f(t)$}
\psfrag{m_f}{$m_f(s,t)$}
\centerline{\includegraphics[height= 3.6 cm]{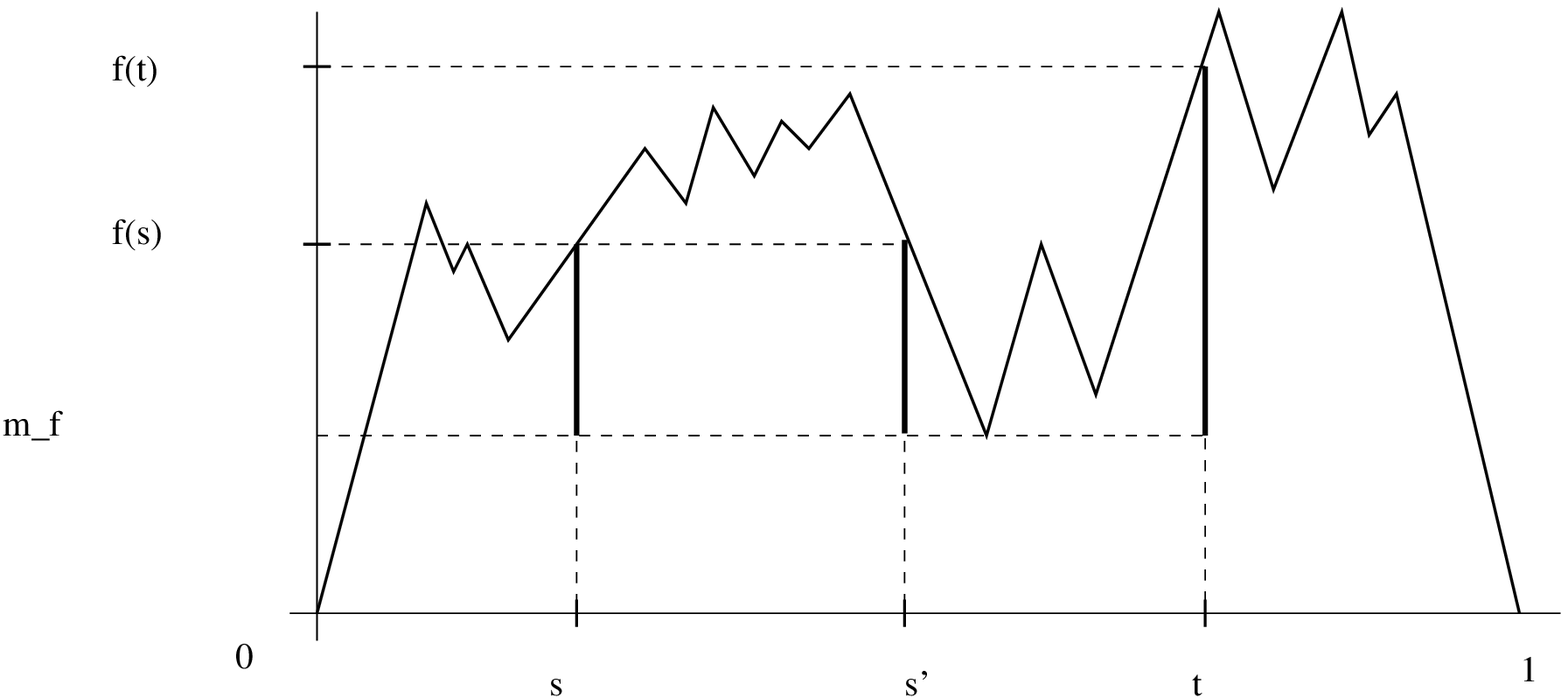}}
\captionn{\label{abst-tree}Graph of a continuous function $f$ satisfying $f(0)=f(1)=0$ and $f(x)\geq 0$ on [0,1]. In this example $s\sous{\sim}f s'$ and the distance $d_f(s,t)=d_f(s',t)=f(s)+f(t)-2\,m_f(s,t)$ is the sum of the lengths of the vertical segments.}
\end{figure}

We adopt the same normalizations as Aldous \cite{ALD, aldous93crt}: the Continuum Random Tree (CRT) ${\cal T}_{2\se}$ 
can be defined as the real tree coded by twice a normalized Brownian excursion $\se = (\se_t)_{t\in[0,1]}$. Indeed, any function $f$ with duration 1 and satisfying moreover $f(0)=f(1)=0$, and $f(x)\geq 0, x\in[0,1]$ may be viewed as coding a continuous tree as follows (illustration can be found on Figure \ref{abst-tree}). For every $s,s'\in[0,1]$, we set
$$m_f(s,s'):=\inf_{s\wedge s'\leq r\leq s\vee s'}f(r).$$
We then define an equivalence relation on $[0,1]$ by setting $s\sous{\sim}f s'$ if and only if $f(s)=f(s')=m_f(s,s')$. Finally we put
\begin{equation}\label{dbe}
d_f(s,s')=f(s)+f(s')-2\,m_f(s,s')
\end{equation}
and note that $d_f(s,s')$ only depends on the equivalence classes of $s$
and $s'$. 

Then the quotient space ${\cal T}_f:=[0,1]/\sous\sim{f}$ equipped with the metric
$d_f$ is a compact $`R$-tree (see e.g. Section 2 of \cite{DuLG}). In other words, it is a compact metric space such that for any two points $\sigma$ and $\sigma'$
there is a unique arc with endpoints $\sigma$ and $\sigma'$ and furthermore
this arc is isometric to a compact interval of the real line.
We view ${\cal T}_f$ as a rooted $`R$-tree, whose root $\rho$ is the equivalence
class of $0$. 

The CRT is the metric space $({\cal T}_{2\se},d_{2\se})$. In addition to the usual genealogical order of the tree, the CRT ${\cal T}_{2\se}$ inherits a lexicographical order (LO) from the coding by $2\se$, in a way analogous to the ordering of (discrete) plane trees from the left to the right. \par
Discrete trees $T$ are now equipped with their graph distances $d_T$.
\begin{pro}\label{crtter}
 The following convergence holds for the GH topology. Under $P^{\ter}_{3n+1}$, 
\[\l(T, \frac{d_T}{\sqrt{3n/2}}\r)\dd ({\cal T}_{2\se},d_{2\se}).\]
\end{pro}
\proof The convergence for the GH topology is a consequence of the convergence for any suitable encoding of trees. The offspring distribution $\nu_\ter:=\frac23\delta_0+\frac13\delta_3$ is critical (in other words has mean 1) and variance $2$. The convergence of rescaled GW trees conditioned by the size proved by Aldous \cite{ALD, aldous93crt}. (See also Le Gall \cite{LGC2} or Marckert \& Mokkadem \cite{MMexc}). ~$\hfill\Box$

\section{Asymptotic behavior of stack-triangulations under $\ut_{2n}$}
\label{res} 

\subsection{Local convergence}

The first aim of this part is to define a map $m_\infty$ built thanks to $\bt_{\infty}^\ter$ with the help of a limiting ``bijection'' analogous to the functions $\Psi_K^\triangle$'s. Some problem arises when one wants to draw or define an infinite map on the plane since we have to deal with accumulation points and possible infinite degree of vertices. We come back on this point later. We now describe a special class of infinite trees -- we call them \sl thin ternary trees \rm -- that will be play an important role further.

\begin{defi}An infinite line of descent in a tree is a sequence $(u_i,i\geq 0)$ such that: $u_0$ is the root $\varnothing$, and $u_i$ is a child of $u_{i-1}$ for any $i\geq 1$. We call thin ternary tree a ternary tree having a unique infinite line of descent $L=(u_i,i\geq 0)$, satisfying moreover $\Gamma(u_n)\sous{\tend}n\infty$ (which will be written $\Gamma(L)=\infty$). The set of thin ternary trees is denoted by $\Tter_\thin$.
\end{defi}

\begin{lem}
\label{supp}
 The support of $P^{\ter}_{\infty}$ is included in $\Tter_\thin$. 
\end{lem}
\proof By construction $L_\infty^\ter$ is an infinite line of descent in $\bt_{\infty}^\ter$ that satisfies clearly a.s. $\Gamma(L_\infty^\ter)=+\infty$. This line is a.s. unique because the sequence $(t(i))$ of grafted trees are critical GW trees and then have a.s. all a finite size.  $~\hfill\Box$\medskip

For any tree $t$, finite or not, denote the $\Gamma-$ball of $t$ of radius $r$ by
\[B_r^\Gamma(t):=\{u~|~ u\in t, \Gamma(u) \leq r \}.\]
\begin{lem}\label{finbal}
 For any tree $t\in \Tter_\thin$ and any $r\geq 0$,  $\#B_r^\Gamma(t)$ is finite. 
\end{lem}
\proof 
Let $L$ be the unique infinite line of descent of $t$. Since $\Gamma(L)=+\infty$, $B_r^\Gamma(t)$ contains only a finite part say $\cro{\varnothing, u}$ of $L$. This part is connected since  $\Gamma$ is non decreasing: if $w=uv$ for two words $u$ and $v$ then $\Gamma(w)\geq \Gamma(u)$. Using again that $\Gamma$ is non decreasing, $B_r^\Gamma(t)$ is contained in $\cro{\varnothing,u}$ union the finite set of finite trees rooted on the neighbors of $\cro{\varnothing,u}$. ~$\hfill\Box$

\begin{pro}\label{geoc} If a sequence of trees $(t_n)$ converges for the local topology to a thin tree $t$, then for any $r\geq 0$ there exists $N_r$ such that for any $n\geq N_r$,
 $B_r^\Gamma(t_n)=B_r^\Gamma(t)$. 
\end{pro}

\proof Suppose that this is not true. Then take the smallest $r$ for which there does not exists such a $N_r$ (then $r\geq 1$ since the property is clearly true for $r=0$).
Let $l_r$ be the length of the longest word in $B_r^\Gamma(t)$. Since $d_L(t_n,t)\leq 1/(l_{r}+1)$ for $n$ say larger than $N'_r$, for those $n$ the words in $t_n$ and $t$ with at most $l_{r}$ letters coincide. This implies that $B_r^\Gamma(t)\subset B_r^\Gamma(t_n)$ and 
that this inclusion is strict for a sub-sequence $(t_{n_k})$ of $(t_n)$. Hence one may find a sequence of words $w_{n_k}$ such that:
$\Gamma(w_{n,k})=r$, $w_{n_k}\in t_{n_k}$, $w_{n_k}\notin t$. Let $w'_{n_k}$ be the smallest (for the LO) elements of $(t_{n_k})$ with this property. In particular, the father $w^f_{n_k}$ of $w'_{n_k}$ satisfies either:\\
$(a)$ $\Gamma(w^f_{n_k})=r-1$ or, \\
$(b)$ $\Gamma(w^f_{n_k})=r$ and then $w^f_{n_k}$ belongs to $B_r^\Gamma(t)$.\\
For $n$ large enough, say larger than $N_{r-1}$, $B_{r-1}^\Gamma(t_n)$ coincides with $B_{r-1}^\Gamma(t)$ (since $r$ is the first number for which this property does not hold). Hence, the set $S_f=\{w^f_{n_k}~|~ n_k \geq N_{r-1}\wedge N'_r\}\subset B_r^\Gamma(t)$ is finite by the previous Lemma. Then the sequence $(w'_{n,k})$ takes its values in the set of children of the nodes of $S_f$, the finite set say $S_r$. Consider an accumulation point $p$ of $(w'_{n_k})$.  The point $p$ is in the finite set $\{w'_{n_k},k\geq 0\}$ and then not in $t$. But $p$ is in $t$ since $t$ contains all (finite) accumulations points of all sequences $(x_n)$, where $x_n\in t_n$. This is a contradiction.  ~$\hfill\Box$

\subsubsection{A notion of infinite map}
\label{nim}
\noindent This section is much inspired by Angel \& Schramm \cite{AS} and Chassaing \& Durhuus \cite[Section 6]{CD}. \medskip

We call infinite map $m$, the embedding of a graph in the plane having the following properties:
\begin{itemize}
\item[$(\alpha)$] it is locally-finite, that is the degree of all nodes is finite,
\item[$(\beta)$] if $(\rho_n, n\geq 1)$ is a sequence of points that belongs to distinct edges of $m$, then accumulation points of $(\rho_n)$ must be outside $m$.
\end{itemize}
This last condition ensures that no face is created artificially. For example, we want to avoid when drawing an infinite graph where each node has degree 2 (an infinite graph line, in some sense) that would create two faces or more: imagine a circular drawing of this graph where the two extremities accumulate on the same point. Avoiding the creation of artificial faces allows to ensure that homeomorphisms of the plane are still the right tools to discriminate similar objects.\par
In the following we define an application $\Psi_{\infty}^\triangle$ that associates with a tree $t$  of $\Tter_\thin$ an infinite map  $\Psi_\infty^\triangle(t)$ of the plane. Before this,  let us make some remarks. Let $t \in \Tter_\thin$, for any $r$,  set  $t(r)$ the tree having as set of internal nodes $B_r^\Gamma=\B_{r}^\Gamma(t)$. We have clearly $d_L(t(r),t)\sous{\to}{r} 0$. Moreover, since $t({r})$ is included in $t({r+1})$, the map $m_{r}=(\Psi^\triangle)^{-1}(t(r))$ is ``included'' in $m_{r+1}$. The quotes are there to recall that we are working on equivalence classes modulo homeomorphisms and that the inclusion is not really defined stricto sensu. In order to have indeed an inclusion, an idea is to use the canonical drawing (see Definition \ref{cd})~: the inclusion ${\cal G}(m_r)\subset {\cal G}(m_{r+1})$ is clear if one uses a history leading to $m_{r+1}$ that passes from $m_r$, which is possible thanks to Proposition \ref{yopp} and the fact that $t(r)\subset t({r+1})$. 
Now $({\cal G}(m_r))$ is a sequence of increasing graphs. Let ${\cal G}_t$ be defined as the map $\cup_r {\cal G}(m_r)$ and having as set of nodes and edges those belonging to at least one of the ${\cal G}(m_r)$.

\begin{pro} 
For any thin tree $t$, the map ${\cal G}_t$ satisfies $(\alpha)$ and $(\beta)$.\end{pro}
\proof 
The first assertion comes from the construction and the finiteness of the balls $B_r^\Gamma$ (by Lemma \ref{finbal}). For the second assertion, just notice that for any $r$, only a unique face of $m_r$ contains an infinite number of faces of ${\cal G}_t$. Indeed, $t(r)$ is included in $t$ and $t$ owns only one infinite line of descent $L$. Hence among the set of fringe subtrees $\{t_u~|~ u\in t(r)\}$ of $t$ (each of them corresponding to the nodes that will be inserted in one of the triangular faces of $m_r$) only one has an infinite cardinality. It remains to check that the edges do not accumulate, and for this, we have only to follow the sequence of triangles $(F_k)$ that contains an infinite number of faces, those corresponding with the nodes of $L$. Moreover, by uniqueness of the infinite line of descent in $t$, the family of triangles $(F_k)$ forms a decreasing sequence for the inclusion. Consider now the subsequence $F_{n_k}$ where $ g(\type(F_{n_k}))=g(\type(F_{n_{k-1}}))+1$. The triangle $F_{n_k}$ has then all its sides different from $F_{n_{k-1}}$. Hence any accumulation points $\rho$ of $(\rho_n)$ (as defined in $(\beta)$) must belong to $\cap F_k$. By the previous argument, $\rho$ does not belong to any side of those triangles, which amounts to saying that $\rho$ lies outside $m$. ~$\hfill\Box$

\begin{pro}\label{contargu}Let $(t_n)$ be a sequence of  trees, $t_n\in\Tter_{3n-2}$, converging for the local topology to a thin tree $t$. Then the sequence of maps $(\Psi_{n}^\triangle)^{-1}(t_n)$ converges to ${\cal G}_t$ for the local topology. 
\end{pro}
\proof If $(t_n)$ converges to $t$ then for any $r$, there exists $n_r$ such that for any $n\geq n_r$, $B_r^\Gamma(t_n)=B_r^\Gamma(t)$. Hence, if $n$ is large enough, $d_L((\Psi_{n}^\triangle)^{-1}(t_{n}),{\cal G}_t)\leq 1/(r+1)$. ~$\hfill\Box$ \medskip

We have till now, work on topological facts, separated in some sense from the probabilistic considerations. It remains to deduce the probabilistic properties of interest.

\subsubsection{A law on the set of infinite stackmaps}

The set $\Tter$ is a Polish space for the topology $d_L$. In such a space, the Skohorod's representation theorem (see e.g.\ \cite[Theorem 4.30]{KAL}) applies: if $(X_n)$ is a sequence of random variables taking their values in a Polish space $S$ and if $X_n\dd X$, then there exists a probability space $\Omega$ where are defined $(\tilde X_n)$ and $\tilde X$ such that, for any $n$, $\tilde X_n\sur{=}{d}X_n$ and $\tilde X\sur{=}{d}X$, and $\tilde X_n\as \tilde X$.
Since $P_{3n-2}^{\ter}$ converges to $P^{\ter}_{\infty}$, there exists a space $\Omega$ on which $\tilde\bt_{n}$ is $P_{3n-2}^{\ter}$ distributed, $\tilde\bt_{\infty}$ is  $P^{\ter}_{\infty}$ and $\tilde\bt_n\as \tilde\bt_{\infty}$. Moreover, thanks to Lemma \ref{supp}, we may assume that a.s. $\tilde\bt_{\infty}$ is a thin tree.

We then work on this space $\Omega$ and use the almost sure properties of $\tilde\bt_{\infty}$. The convergence in distribution  of our theorem will be a consequence of the a.s. sure convergence on $\Omega$.
\begin{defi}
We denote by $`P^{`3}_{\infty}$ the distribution of $m_\infty:={\cal G}_{\bt_{\infty}}$.
\end{defi}
A simple consequence of Proposition \ref{contargu} is the following assertion. Since $d_L(\tilde\bt_n, \tilde\bt_{\infty})\as 0$ then 
\begin{equation}\label{asconvcarte}
d_L\l((\Psi_{n}^\triangle)^{-1}(\tilde\bt_{n}), {\cal G}_{\tilde\bt_{\infty}}\r)\as0.
\end{equation}
 This obviously implies the following result.
\begin{theo}\label{loctri} $(\ut_{2n})$ converges weakly to $`P^{`3}_{\infty}$ for the topology of local convergence. 
\end{theo}

\subsection{Asymptotic under the Gromov-Hausdorff topology}
\label{asGH}
We begin with a simple asymptotic result concerning the function $\Gamma$ defined in Section \ref{fond-bij}. 
\begin{lem}\label{cou}
Let $(X_i)_{i\geq 1}$ be a sequence of random variables uniform in $\Sigma_3=\{1,2,3\}$, and independent. Let $W_n$ be the word $(X_1,\dots,X_n)$.\\
$(i)$ $n^{-1}{\Gamma(W_n)}\as\Gamma_\triangle$ 
where 
\begin{equation}
\Gamma_\triangle:=2/11.
\end{equation}
$(ii)$ $`P( |\Gamma(W_n)-n\Gamma_\triangle|\geq n^{1/2+u})\sous{\to}n 0$ for any $u>0$.
\end{lem}
\proof If $W$ is the infinite sequence $(X_i)$, clearly $\tau_2(W)\sim \Geo(1/3)$ and for $i\geq 3$, the $(\tau_i(W)-\tau_{i-1}(W))'s$ are i.i.d., independent also from $\tau_2$, and are distributed as $1+G_1+G_2$ where $G_1\sim \Geo(1/3)$ and $G_2\sim \Geo(2/3)$ [the distribution $\Geo(p)$ is $\sum_{k\geq 1}p(1-p)^{k-1}\delta_k$]. It follows that $`E(\tau_i(W)-\tau_{i-1}(W))=11/2$ for $i\geq 3$ and $`E(\tau_2(W))=3<+\infty$. By the renewal theorem assertion $(i)$ holds true. For the second assertion, 
write
\[\{|\Gamma(W_n)-n\Gamma_\triangle|\geq n^{1/2+u}\}=\{\tau_1+\dots+\tau_{n\Gamma_\triangle+n^{1/2+u}}\leq n\}\cup\{\tau_1+\dots+\tau_{n\Gamma_\triangle-n^{1/2+u}}\geq n\}.\]
By the Bienaymé-Tchebichev inequality the probability of the events in the right hand side goes to 0. $\hfill\Box$\medskip

For every integer $n \geq 2$, let $M_n$ be a random rooted map under $\ut_{2n}$. Denote by $m_n$ the set of vertices of $M_n$ and by $D_{m_n}$ the graph distance on $m_n$. We view $(m_n, D_{m_n})$ as a random variable taking its values in the space of isometric classes of compact metric spaces.   
\begin{theo}\label{youp}
Under $\ut_{2n}$,   
\[\l(m_n,\frac{D_{m_n}}{\Gamma_\triangle\sqrt{3n/2}}\r)\dd ({\cal T}_{2\se}, d_{2\se}),\] 
for the Gromov-Hausdorff topology on compact metric spaces.
\end{theo}

This theorem is a corollary of the following stronger Theorem stating the convergence of maps seen as parametrized metric spaces.
In order to state this theorem, we need to parametrize the map $M_n$. 
The set of internal nodes of $m_n$ inherits of an order, the LO on trees, thanks to the function $\Psi^\triangle_n$. Let $u(r)$ be the $r$th internal node of $m_n$ for $r\in\{0,\dots,n-1\}$. Denote by $d_{m_n}(k,j)$ the distance between $u(k)$ and $u(j)$ in $m_n$. We need in the following theorem to interpolate $d_{m_n}$ between the integer points to obtain a continuous function. Any smooth enough interpolation is suitable. [For example, define $d_{m_n}$ as the plane interpolation on the triangles with integer coordinates of the form $(a,b),(a+1,b),(a,b+1)$ and $(a,b+1),(a+1,b+1),(a+1,b)$].
\begin{theo}\label{paramversion}
Under $\ut_{2n}$,  
\begin{equation}
\label{par}
\l(\frac{d_{m_n}(ns,nt)}{\Gamma_\triangle\sqrt{3n/2}}\r)_{(s,t)\in[0,1]^2}\dd \l(d_{2\se}(s,t)\r)_{(s,t)\in[0,1]^2},
\end{equation}
where the convergence holds in $C[0,1]^2$ (even if not indicated, the space $C[0,1]$ and $C[0,1]^2$ are equipped with the topology of uniform convergence).
\end{theo}
The proof of this Theorem is postponed to Section \ref{rp}.\par
The profile $\Prof_m:=(\Prof_m(t),t\geq 0)$ of a map $m$ with root vertex $E_0$ is the càdlàg-process
\[\Prof_m(t)=\#\{u \in V(m)~|~, d_m(E_0,u)\leq t\}, \textrm{ for any }t\geq 0.\] 
The radius $R(m)=\max\{d_m(u,E_0)~|~u \in V(m)\}$ is the largest distance to the root vertex in $m$.\par
As a corollary of Theorem \ref{youp} or Theorem \ref{paramversion}, we have:
\begin{cor}
Under $\ut_{2n}$, the process 
\begin{equation}\label{pozer}
\l(n^{-1}\Prof_{m_n}(\Gamma_\triangle\sqrt{3n/2}\,v)\r)_{v\geq 0}\dd \l(\int_0^v l^x_{2\se}dx\r)_{v\geq 0}
\end{equation}
where $l^x_{2\se}$ stands for the local time of twice the Brownian excursion $2\se$ at position $x$ at time 1, and where the convergence holds in distribution in the set $D[0,+\infty)$ of  c\`adl\`ag functions endowed with the Skohorod topology. Moreover 
\[\frac{R(m_n)}{\Gamma_\triangle\sqrt{3n/2}}\dd 2\max \se\] 
\end{cor}
\proof Let $D_n(s)=\frac{d_{m_n}(ns,0)}{\Gamma_\triangle\sqrt{3n/2}}$ be the interpolated distance to $E_0$. 
By \eref{par}, $({D_n(s)})_{s\in[0,1]}\dd(2\se(s))_{s\in[0,1]}$ in $C[0,1]$.
By the Skohorod's representation theorem (see e.g.\ \cite[Theorem 4.30]{KAL}) there exists a space $\Omega$ where a copy $\tilde{D}_n$ of $D_n$, and a copy $\tilde{\se}$ of $\se$ satisfies $\tilde{D}_n\as 2\tilde{\se}$ in $C[0,1]$. We work from now on on this space, and write $\widetilde{\Prof}_n$ the profile corresponding to $\tilde{D}_n$. For any $v$ such that $\Gamma_\triangle\sqrt{3n/2}\,v$ is an integer,
\[n^{-1}\widetilde{\Prof}_{n}(\Gamma_\triangle\sqrt{3n/2}\,v) = \int_0^1
\1_{\tilde D_{n}(s)\leq v}\,ds.\]
For every $v$, a.s., $\int_0^1 \1_{\tilde{D}_{n}(s)\leq v}\,ds \to \int_0^v l^x_{2\tilde{\se}}\,dx$. To see this, take any $`e>0$ and check that $\|\tilde{D}_n-2\tilde{\se}\|_{\infty}\to 0$ yields 
\begin{equation}\label{era}
\int_0^1 \1_{2\tilde{\se}(s)\leq v-`e}\,ds \leq \int_0^1
\1_{\tilde{D}_{n}(s)\leq v}\,ds\leq \int_0^1 \1_{2\tilde{\se}(s)\leq v+`e}\,ds.
\end{equation}
Since the Borelian measure $\mu_{2\se}(B)=\int_0^1 \1_{2e(s)\in B}\,ds$ has no
atom a.s., $v\to \int_0^v l^x_{2\se}dx$ is continuous and non-decreasing.
Hence since  $v\to\int_0^1 \1_{D_{n}(s)\leq v}\,ds$ is non decreasing and by
\eref{era} we have  $\int_0^1 \1_{D_{n}(s)\leq v}\,ds \to \int_0^v l^x_{2\se}\,dx$ a.s. for any $v\geq 0$. Thus, $(v\to \int_0^1 \1_{D_{ns}\leq v}\,ds) \to (v\to\int_0^v l^x_{2\se}\,dx)$ in $C[0,1]$. This yields the convergence of $\Prof_{m_n}$ as asserted in \eref{pozer}. \par
For the second assertion, note that $f\to \max f$ is continuous on $C[0,1]$. Since $\tilde{D}_n\as 2\tilde{\se}$ then $\max \tilde{D}_n\as \max 2\se$, and then also in distribution. ~$\hfill\Box$

\subsection{Asymptotic behavior of the typical degree}

\begin{pro}\label{prop:degre}
Let $m_n$ be a map $\ut_{2n}$ distributed, $u(1)$ the first node 
inserted in $m_n$, and  ${\bf u}$ be a random node chosen uniformly among 
the internal nodes of $m_n$.\\
$(i)$ $\deg_{m_n}(u(1))\dd X$ where for any $k\geq 0$, $`P(X=k+3)=\frac{k}{k+3}
\binom{2k+2}{k}
\frac{2^{k+3}}{3^{2k+3}}
$ .\\
$(ii)$ $\deg_{m_n}({\bf u})\dd Y$ where for any $k\geq 0$, 
$`P(Y=k+3)=\frac{1}{k+3}\binom{2k+2}{k}\frac{2^{k+3}}{3^{2k+2}}
$.
\end{pro}

\begin{lem}\label{lem:degre}Let $T$ be a random tree under $\uter_{3n+1}$ and ${\bf u}$ be chosen uniformly in $T^\circ$. We have
$|T_{\bf u}|\dd {\bf K}$ where
$`P({\bf K}=3k+1)=\frac{2^{2k+1}}{3^{3k}(3k+1)}\binom{3k+1}{k}, \textrm{ for }k\geq 1. $
Moreover, conditionally on $|T_{\bf u}|=m$, $T_{{\bf u}}$ has the uniform distribution in $\Tter_m$.
\end{lem}
\proof Consider 
\be \Tters_{3n+1}:=\{(t,u)~ |~ t\in\Tter_{3n+1}, u\in t^\circ\},~~~
\Tterb_{3n+1}:=\{(t,u)~ |~ t\in\Tter_{3n+1}, u\in \partial t\}
\ee
the set of ternary trees with a distinguished internal node, resp. leaf. For any tree $t$ and  $u\in t$ set $t[u]=\{v
\in t~|~ v \textrm{ is not a descendant of }u\}$.  Each element
$(t,u)$ of $\Tters_{3n+1}$ can be decomposed bijectively as a pair $[(t[u],u),t_u]$ where $(t[u],u)$ is a tree with a marked leaf, and $t_u$ is a ternary tree having at least one internal node. 
Hence, For any $n$, the function $\rho$ defined by $\rho(t,u):=[(t[u],u),t_u]$ is a
bijection from  $\Tters_{3n+1}$ onto $\bigcup_{k=1}^{n}\l(\Tterb_{3(n-k)+1}\times \Tter_{3k+1}\r)$. \par
Since the trees in $\Tter_{3n+1}$ have the same number of internal nodes,
choosing a tree $T$ uniformly in $\Tter_{3n+1}$ and then a
node $u$ uniformly in $T^\circ$, amounts to choosing a
marked tree $(T, u)$ uniformly in $\Tters_{3n+1}$. We then have, for
any fixed $k$,
\begin{eqnarray}
  \uter_{3n+1}(|T_{\bf u}|=3k+1)&=&\#\Tterb_{3(n-k)+1}
  \#\Tter_{3k+1} \l(\#\Tters_{3n+1}\r)^{-1}.
  \label{eq:subtree}
\end{eqnarray}
When $n\to+\infty$, this tends to the result announced in the Lemma, using  $\#\Tterb_{3m+1}=(2m+1)\#\Tter_{3m+1}$ and $\#\Tters_{3n+1}=n\#\Tter_{3n+1}$ and 
\begin{equation}\label{enu-ter}
\#\Tter_{3n+1}=\frac{1}{3n+1}\binom{3n+1}{n}\sim \sqrt\frac{3}{\pi}\frac{3^{3n}}{2^{2n+2}n^{3/2}}.
\end{equation}
Since $\sum_{k\geq 1} \frac{2^{2k+1}}{3^{3k}(3k+1)}\binom{3k+1}{k}=1$, we have
indeed a convergence in distribution of $\deg_T({\bf u})$ under $\uter_{3n+1}$
to ${\bf K}$.  The second assertion of the Lemma is clear. ~$\Box$ \medskip

\begin{figure}[htbp]
\centerline{\includegraphics[height= 2.5 cm]{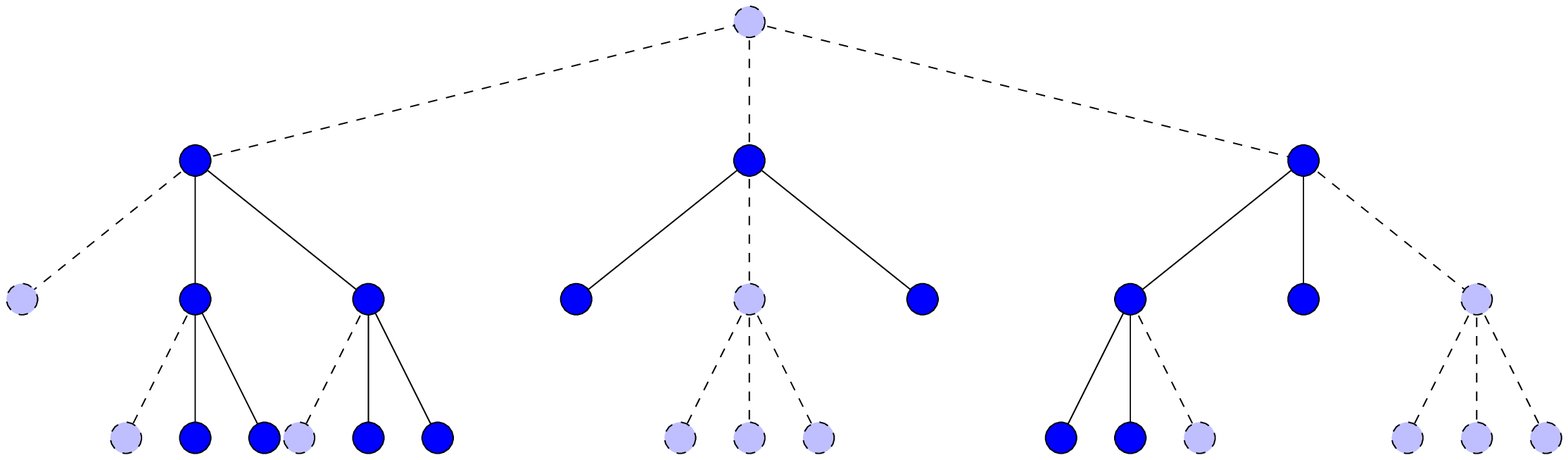}}
\captionn{\label{fig:dec} A ternary tree $t$ and $t^{deg}$. Plain 
vertices belong to $t^{deg}$.} 
\end{figure}

\prooff{Proposition \ref{prop:degre}}
As illustrated on Figure \ref{fig:dec}, for any  $t\in\Tter$, we let
$$
t^{deg}:= \{v~|~ v\in t, v \in 1L^\star_{2,3}\cup 2L^\star_{1,3}\cup 3L^\star_{1,2}\}. 
$$
In general $t^{deg}$ is a forest of three pseudo-trees: pseudo here means that the connected components of $t^{deg}$ have a tree structure but do not satisfies the first and third points in Definition \ref{tree}. For sake of compactness, we will however up to a slight abuse of language call these three pseudo-trees, binary trees (combinatorially their are binary trees).
 
\noindent $(i)$ Let $T$ be a tree $\uter_{3n-2}$ distributed and $m=(\Psi^\triangle_n)^{-1}(T)$. By Proposition \ref{yop}, 
\begin{equation}
\deg_m(u(1))=3+\#(T^{deg}\cap T^\circ),
\end{equation}
or in other words $\ut_{2n}(\deg(u(1)=k)=\uter_{3n-2}(|T^{deg}|=2k+3)$.
Each ternary tree $t$ not reduced to the root vertex can be decomposed in a unique way as a pair $(t^{deg},f)$ where  $f:=(t(1), \dots, t(k))\in (\Tter)^k$ is a forest of ternary trees, and $k=\#(t^{deg}\cap t^\circ)$. Let $\Fbin^{n}(k)$ 
(resp $\Fter^n(k)$) be the set of forests composed with $n$ binary (resp. ternary) trees and total number of nodes $k$. For $0\leq k <n-1$, we get:
\begin{equation}
\uter_{3n-2}(|T^{deg}|=2k+3)=\frac{\#\Fbin ^3 (2k+3) 
\#\Fter^k(3n-2k-6)}
{\#\Tter_{3n-2}}.
\end{equation}
A well known consequence of the rotation/conjugation principle is that 
\begin{equation}
\# \Fbin ^m (n)  = \frac{m}{n}\binom{n}{(n+m)/2}, ~~~\textrm{ and }\# \Fter ^m (n)  = \frac{m}{n}\binom{n}{(n-m)/3}
\end{equation}
with the convention that $\binom{a}{b}$ is 0 if $b$ is negative or non integer.
We then have
\begin{equation}\label{eq:deg_root}
\uter_{3n-2}(|T^{deg}|=2k+3)=
\frac{	\frac{3}{2k+3}
\binom{2k+3}{k}
\frac{k}{3n-2k-6}
\binom{3n-2k-6}{n-k-2}}
{\frac{1}{3n-2}
\binom{3n-2}{n-1}}.
\end{equation}
We get
$\ut_{2n}(\deg_{m_n}(u(1))=k+3)
\sous{\to}{n}
\frac{k}{k+3}
\binom{2k+2}{k}
\frac{2^{k+3}}{3^{2k+3}}$, limit which is indeed a probability distribution.\\
$(ii)$ Now let $m_n$ be $\ut _{2n}$ distributed and $\bf{u}$ be a  uniform internal node of $m_n$. Let $T=\Psi^\triangle(m_n)$ and $\mathbf{u}'$ be the internal node of $T$ corresponding to $\mathbf{u}$.  We have this time $\ut_{2n}(\deg_m({\bf u})=k)=\uter_{3n-2}(|T^{deg}_{\bf u'}|=2k+3)$.
First by a simple counting argument, 
\[
\uter_{3n-2}\l(|T^{deg}_{\mathbf{u}'}|=2k+3~ \big|~ |T_{\bf{u'}}| = 3j-2\r)=
\uter_{3j-2}\l(|T^{deg}|=2k+3\r). 
\]
Conditioning on $|T_{\bf{u'}}|$, using Formulas \eref{eq:subtree} and (\ref{eq:deg_root}) we get after simplification
\begin{equation}\label{eq:degre}
  \uter _{3n-2}\l(|T_{\mathbf{u}'}^{deg}|=2k+3\r)=  \sum _{j \geq k+2}  q_{n,k,j}\end{equation}
where 
\[q_{n,k,j}=  \l(\1_{j\leq n}\r)~\frac{3}{2k+3}\binom{2k+3}{k}\frac{\frac{k}{3j-2k-6}\binom{3j-2k-6}{j-k-2}
 \binom{3(n-j)}{n-j}}{\binom{3n-3}{2n-1}}.\]
We have $\lim_n\binom{3(n-j)}{n-j}/
  {\binom{3n-3}{2n-1}}=2^{2j-1}/3^{3j-3}$, and thus 
\begin{equation}\label{v_k}
  \sum _{j \geq k+2} ^{\infty}\lim_n q_{n,k,j}=
    \frac{3}{2k+3}\binom{2k+3}{k}\frac{2^{k+3}}{3^{2k+3}},
 \end{equation}
which is the probability distribution announced. To end the proof we have to explain why the exchange $\lim_n$ and $\sum_{j\geq k+2}$ is legal. Recall the Fatou's Lemma: if $(f_i)$ is a sequence of non-negative measurable functions, then $\int \liminf_{n\to\infty} f_n\,d\mu \le \liminf_{n\to\infty} \int_S f_n\,d\mu.$
Set 
\[v_k:= \sum_{j \geq k+2} ^{\infty} \liminf_n q_{n,k,j},~~~~ u_k  :=\liminf_n  \sum _{j \geq k+2}  q_{n,k,j}.\]
The sequence $(v_k)$ has been computed in \eref{v_k} and sum to $1$. 
By Fatou's Lemma, $u_k\geq v_k.$ By Fatou's Lemma again,
\[\sum_k u_k \leq \liminf_n \sum_k \sum _{j \geq k+2} q_{n,k,j}=1.\] We deduce that $u_k=v_k$.~$\Box$

\section{Asymptotic behavior of stack-triangulations under $\qt_{2n}$}
\label{res2}

We first work on ternary trees under  $\qter_{3n-2}$ (recall the content of Section \ref{ind-dist}).
\begin{pro}\label{dic}
Let $\bt$ be a random tree under $\qter_{3n-2}$, and ${\bf u}$ and ${\bf v}$ be two i.i.d. random variables uniform in $\bt^\circ$, the set of internal nodes of $\bt$. Let ${\bf w}={\bf u}\wedge {\bf v}$.\\
1) We have
$\l(\frac32\log n\r)^{-1/2}\l(|{\bf u}|-\frac{3}2\log n,|{\bf v}|-\frac{3}2\log n\r)\dd (N_1,N_2)$
where $N_1$ and $N_2$ are independent centered Gaussian r.v. with variance 1.\\
2) Let ${\bf a}, {\bf b}\in \{1,2,3\}$, with ${\bf a}\neq {\bf b}$ and  ${\bf u}^\star, {\bf v}^\star$ the (unique) words such that
\[{\bf u}={\bf wau}^\star \textrm{ and } {\bf v}={\bf wbv}^\star.\]
Conditionally to $(|{\bf u}^\star|,|{\bf v}^\star|)$ (their lengths) ${\bf u}^\star$ and ${\bf v}^\star$ are independent random words composed with $|{\bf u}^\star|$ and $|{\bf v}^\star|$ independent letters uniformly distributed in $\Sigma_3=\{1,2,3\}$.
\end{pro}

This Proposition is more or less part of the folklore. In Bergeron \& al \cite{BFS}, in particular in Theorem 8 and Example 1 p.7, it is proved that
\begin{equation}\label{re}
\l(\frac32\log n\r)^{-1/2}\l(|{\bf u}|-\frac{3}2\log n\r)\dd N_1.
\end{equation}
The fact that $|{\bf u}^\star|$ and $|{\bf v}^\star|$ behave as $|{\bf u}|$ and are asymptotically independent comes from that ${\bf w}$ is close to the root, and also from the linear size of the two subtrees rooted in ${\bf w}$ containing ${\bf u}$ and ${\bf v}$ (the normalizations in Formula \eref{re} are asymptotically insensible to the use of  $an$ instead of $n$), and are, given their size, increasing trees with these sizes. The uniformity of the letters comes from a symmetry argument. 
Below we present a formal proof of this proposition using a ``Poisson-Dirichlet fragmentation'' point of view, very close to that used in Broutin \& al. \cite[Section 7]{BDMD} where the height of increasing trees is investigated. 
They show that 
in increasing trees 
the asymptotic proportion $n^{-1}(|t_1|,\dots,|t_d|)$  of nodes in the subtrees of the root are given by a Poisson-Dirichlet distribution. The point of view developed below is slightly different, since we first take a Poisson-Dirichlet fragmentation and then show that the fragmentation tree is distributed as an increasing tree, leading then at once to the convergence of $n^{-1}(|t_1|,\dots,|t_d|)$. The following Subsection is mostly contained in the more general work of Dong \& al. \cite{DGM} (particularly Section 5). We give a straight exposition below for the reader convenience, in a quite different vocabulary.

\subsection{Poisson-Dirichlet fragmentation}

We construct here a representation of the distribution $\qt_{3K-2}$ as the distribution of the underlying tree of a fragmentation tree. Let begin with the description of the deterministic fragmentation tree associated with a sequence of choices ${\bf b}=(b_i)_{i\geq 1}$, $b_i\in[0,1]$ and a sequence ${\bf y}=(y^u)_{u\in {\cal W}_3}$ (indexed by the infinite complete ternary tree), where for each $u$,
\[y^u=(y^u_1,y^u_2,y^u_3)\]
where for any  $i\in\{1,2,3\}$ and $u\in {\cal W}_3$, $y^u_i>0$ and $\sum_{i=1}^3 y^u_i=1$. The sequence $(y^u)$ may be thought as the fragmentation structure associated with the tree.

With these two sequences we associate a sequence $F_n=F(n,{\bf b},{\bf y})$ of ternary trees with $3n+1$ leaves, where each node is marked with an interval as follows.\\
-- At time 0, $F_0$ is the tree $\{\varnothing\}$ (reduced to the root) marked by $I_{\varnothing}=[0,1)$.\\
-- Assume now  that $F_i$ is built, and is a ternary tree with $3i+1$ nodes each marked with an interval included in $[0,1)$, and such that the leaves-intervals $(I_u,u\in \partial T_i)$ form a partition of $[0,1)$. Then the tree $F_{i+1}$ is obtained from $F_i$ as follows. Consider $u^\star$ the leaf whose associated interval $I_{u^\star}$ contains $b_{i+1}$. Give to $u^\star$ the 3 children $u^\star1,u^\star2,u^\star 3$. Now split the interval $I_{u^\star}$ into $(I_{u^\star1},I_{u^\star2}, I_{u^\star3})$ with respective size proportions given by $y^{u^\star}$: if $I_{u^\star}=[a,b)$ then set  $I_{u^\star i}=[a+(b-a)\sum_{j=1}^{i-1}y^{u^\star}_j,a+(b-a)\sum_{j=1}^{i}y ^{u^\star}_j)$ for every $i\in\{1,2,3\}$. 
Let $\Omega_{\cal F}$ be the set of fragmentation trees (a tree where each node is marked by an interval).
We define the application $\pi$ from $\Omega_{\cal F}$ to $\Tter$ the application sending a fragmentation tree $F$ to its underlying tree $\pi(F)$, that is the tree $F$ without marks.\par
We now let ${\bf b}$ and ${\bf y}$ be random. 
For $d\geq 2$ consider the simplex
\[\Delta_{d-1}=\l\{x=(x_1,\dots,x_{d})~|~ x_i\geq 0 \textrm{ for every }i\in\{1,\dots,d\} \textrm{ and } \sum_{i=1}^{d}x_i=1\r\}.\] The $d-1$-dimensional Dirichlet distribution with parameter $\alpha\in(0,+\infty)$, denoted $\Dir_{d-1}(\alpha)$, is the probability measure (on $\Delta_{d-1}$) with density 
\begin{equation}
\mu_{d,\alpha}(x_1,\dots,x_{d}):=\frac{\Gamma(d\alpha)}{\Gamma(\alpha)^{d}}\,x_1^{\alpha-1}\dots x_{d}^{\alpha-1}
\end{equation}
with respect to $dS_d$ the uniform measure on $\Delta_{d-1}$.
Consider the following discrete time process $({\bf F}_n)$ where ${\bf F}_n=F(n,{\bf B},{\bf Y})$,  ${\bf B}$ is a sequence of i.i.d. random variables uniform on $[0,1]$, and ${\bf Y}=(Y^u)_{u\in {\cal W}_d}$ is a sequence of i.i.d. r.v. with $\Dir_{d-1}(\alpha)$ distribution (independent from ${\bf B}$). When like here, the choice of the interval that will be fragmented is equal to the size of the fragment, the fragmentation is said to be biased by the size. 
\begin{pro}If $d=3$ and $\alpha=\frac1{d-1}$ for any $K\geq 1$ the distribution of $\pi({\bf F}_K)$ is $\qter_{3K-2}$.
\end{pro} 
\begin{note}For any $d\geq 2$, the distribution of the underlying fragmentation tree is a distribution on $d$-ary tree similar to $\qter_{3K-2}$: it corresponds to $d$-ary increasing trees, and can also be constructed thanks to successive insertions of internal nodes uniformly on the existing leaves.
\end{note}

\proof Let ${\bf t}^{(K)}=\pi({\bf F}_K)$. Due to the recursive structure of fragmentation trees, the distribution of the size of the subtrees $(|{\bf t}^{(j)}_1|,|{\bf t}^{(j)}_2|,|{\bf t}^{(j)}_3|)$ for every $j\leq K$, characterizes the distribution of ${\bf t}^{(K)}$. Knowing $Y^\varnothing=(Y^\varnothing_1,Y^\varnothing_2,Y^\varnothing_3)$, the distribution of $(|({\bf t}^{(K)}_1)^\circ|,|({\bf t}^{(K)}_2)^\circ|,|({\bf t}^{(K)}_3)^\circ|)$
is multinomial $(K-1,Y^\varnothing_1,Y^\varnothing_2,Y^\varnothing_3)$; indeed, insertions are ruled out by the number of variables $(B_i,i\leq K-1)$ belonging to each of the intervals $I_{i}=[\sum_{j=1}^{i-1}Y^\varnothing_j, \sum_{j=1}^{i}Y^\varnothing_j)$ for any $i\in \{1,2,3\}$.\par
Let us integrate this. We have
\begin{equation}
\label{prem}
P\l(|t_{i}^\circ{}^{(K)}|=k_i, i\in\{1,2,3\}\r)=\int_{\Delta_{2}}\binom{K-1}{k_1,k_2,k_3}x_1^{k_1}x_2^{k_2} x_3^{k_3}\mu_{3,\frac{1}{2}}(x_1,x_2,x_{3})dS_3(x_1,x_2,x_3)
\end{equation}
for any non negative integers $k_1,k_2,k_3$ summing to $K-1$. This leads to
\ben\label{momentdir}
P\l(|t_{i}^\circ{}^{(K)}|=k_i, i\in\{1,2,3\}\r)=\binom{K-1}{k_1,k_2,k_3}\frac{\Gamma(3/2)}{\Gamma(1/2)^{3}}\frac{\prod_{i=1}^3 \Gamma(k_i+1/2)}{\Gamma(k_1+k_2+k_3+3/2)}\een 
The comparison with $\qter$ is done as follows. Let count the number of constructions leading to a tree $t$ such that $|t_{i}^\circ|=k_i, i\in\{1,2,3\}$. The sum of the number of histories of the trees with $m$ internal nodes is $N_m:=\prod_{i=0}^{m-1} (2i+1)$ since each time the number of leaves increases by 2. 
Hence 
\[\qter_{3K-2}(|{t}_{i}^{\circ}|=k_i,i\in\{1,2,3\})=\binom{K-1}{k_1,k_2,k_3}\frac{\prod_{i=1}^3 N_{k_i}}{N_K}.\]
A simple computation shows that this is proportional to \eref{momentdir}. Since two proportional distribution are equals, we have the result.~$\hfill\Box$ \medskip

\prooff{ Proposition \ref{dic}} \rm In a size biased fragmentation process where the fragmentation measure does not charge 0, the maximal size of the fragments goes a.s. to 0 when the time goes to $+\infty$. 
Hence for any $`e>0$ and $`e'>0$ fixed, for $r$ large enough,
\[P(\max\{|I_u|,u\in \partial \pi(F_r)\}\leq `e)\geq 1-`e'.\] 
Now let us work conditionally on  $E:=\{(|I_u|\leq `e, u\in \partial\pi(F_r))\}$, the event that all fragments have size smaller than $`e$ at time $r$, and consider the fragmentation tree $\bt^{(n)}:=\pi(F_n)$ at time $n$, for $n\geq r$. The vector $(|(\bt^{(n)}_u)^\circ|,u\in \partial \bt^{(r)})$ [giving the number of internal nodes in the fringe subtrees at time $n$] is multinomial $(n-r,(|I_u|, u\in \partial F_r))$. Hence conditionally on ${\bf u},{\bf v} \notin \bt^{(r)}$ (which happens with probability $(n-r)^2/n^2\geq 1-`e$ when $n$ is large), the probability that $\bf u$ and $\bf v$ are chosen in the same subtree is given by $\sum_{u\in\partial F_r} |I_u|^2\leq \max |I_u|\sum_{u\in\partial F_r}|I_u|=\max |I_u|\leq `e$. 
In this case, the height $|{\bf w}|$ (where ${\bf w}={\bf u}\wedge{\bf v}$) is smaller than  $r$ (since the height of $\bt^{(r)}$ is smaller than $r$). It remains to say that conditionally on $({\bf w}, I_{\bf w}, y^{\bf w})$, the strong law of large numbers ensures that the subtrees $(\bt^{(n)}_{{\bf w}i},i=1,2,3)$ (those rooted at the children of ${\bf w}$),  have an asymptotic linear size with $n$, when $n$ goes to $+\infty$  (since the number of $B_i$'s, $r< i\leq n$ fallen in a given interval follows a binomial distribution). Moreover conditionally on their sizes, they are copies of fragmentations trees and then behaves, in terms of shape, as increasing trees. Moreover, since ${\bf u}$ and ${\bf v}$ are chosen uniformly in $\bt^\circ$, knowing that ${\bf u}$ (and ${\bf v}$) is in a given subtree, yields that it is uniformly distributed in this subtree. Then Formula \eref{re} applies. This allows to get $(1)$; then $(2)$ follows by a symmetry argument.
~$\hfill\Box$\medskip

The following theorem may be considered as the strongest result of this section. 
\begin{theo}\label{metconv} Let $M_n$ be a stack-triangulation under $\qt_{2n}$. Let $k\in \mathbb{N}$ and ${\bf v}_1,\dots,{\bf v}_k$ be $k$ nodes of $M_n$ chosen independently and uniformly among the internal nodes of $M_n$. We have 
\[\l(\frac{D_{M_n}({\bf v}_i,{\bf v}_j)}{3\Gamma_\triangle\log n}\r)_{(i,j)\in\{1,\dots,k\}^2}\proba \l(1_{i\neq j}\r)_{(i,j)\in\{1,\dots,k\}^2}\] the matrix of the discrete distance on a set of $k$ points.
\end{theo}
This is consistent with the computations of Zhou \cite{ZYW} and   Zhang \& al \cite{ZRC}.

\proof This is a consequence of Lemma \ref{cou} and the pairwise convergence provided by Proposition \ref{dic} (asymptotically the distance in the tree between two random nodes ${\bf u}$ and ${\bf v}$ is asymp. around $3\log n$, and the letters of ${\bf u}^\star$ and ${\bf v}^\star$ are independent) together with Lemma \ref{cou}.~$\hfill\Box$

We give now some indications about the limiting behavior of triangulations under the law $\qt_{2n}$.

\subsection{Some features of large maps under $\qt$}
\label{azd}
 Some asymptotic results allowing to understand the behavior of large maps under $\qt$ can also be proved using the fragmentations processes. In particular using that the size of a subtree rooted on a given node $u$ evolves (asymptotically) linearly in time (this is due, as said before, to the rate of insertions of nodes in $T_u$ which is constant and given by $|I_u|$), the same results holds true for a fixed face in the triangulation. Moreover, the length $|I_u|$ is the product of $|u|$ marginals of Poisson-Dirichlet random variables. Hence $N_n(f)$ the number of internal nodes present in the canonical face $f$ at time $n$ behaves as follows: $n^{-1}N_n(f)$ converges a.s. toward a random variable $N_f$ almost surely in $(0,1)$. This fragmentation point of view allows to prove much more as the a.s. joint convergence  of $n^{-1}(N_{f_1},\dots,N_{f_{k}})$ for the (disjoint or not) faces $f_i$ of $m_j$ toward a limiting random variable taking its value in $`R^k$, and whose limiting distribution may be described in terms of product of Poisson-Dirichlet random variables. \par

The degree of a node may also be followed when $n$ goes to $+\infty$. If $v(j)$ denotes the $j$th node inserted in $m_n$, one may prove that $\deg(v(j))$ goes to infinity with $n$. The degree of a node follows indeed a simple Markov chain since it increases if and only if a node is inserted in a face adjacent to $v(j)$ and this occurs with a probability equals to $\deg(v(j))$ divided by the current number of internal faces. Denoting by $D_j^n$ the degree of $\deg(v(j))$ at time $n$ (recall that $D_j^j=3$),  under $\qt_{2n}$, we have that 
for $n> j$ and $k\geq 3$, conditionally on $D_j^n$
\begin{equation}\label{evold}
D_{j}^{n+1}=D_j^{n}+B\l(D_j^{n}/(2n-1)\r)
\end{equation}
where we have denoted by $B(p)$ a Bernoulli random variable with parameter $p$ (in other words $\qt_{2(n+1)}(\deg(v(j))=k+1)=\frac{k}{2n-1}\qt_{2n}(deg(v(j))=k)+\frac{2n-k-2}{2n-1}\qt_{2n}(deg(v(j))=k+1)).$

This chain has the same dynamics as the following simple model of urn. Consider an urn with 3
white balls and $2j-2$ black balls at time 0. At each step pick a
ball and replace it in the urn. If the picked ball is white then add one white ball and one black ball, and
if it is black, add two black balls. The number $N_j^t$ of white balls at time $t$ has the same law as $D_j^{j+t}$ (the number of black balls behaves as the number of finite faces of $m_{j+t}$ not incident to $v(j)$). 
This model of urn has been studied in Flajolet \& al. \cite[p.94]{FDP} (to use their results, take $a_0=3$, $b_0=2j-2$, $\sigma=2$, $\alpha=1$ and replace $n$ by $n-j$). For example, we derive easily from their results the following proposition.
\begin{pro}
Let $m_n$ be a map $\qt_{2n}$ distributed and $v(j)$ the j-th node inserted, for
$n>j$ and $1\leq k \leq n-j$, we get
\[
\qt_{2n}(deg_{m_n}(v(j))=k+3)=\frac{\Gamma(n-j+1)\Gamma(j+\frac12)}{\Gamma(n+\frac12)}\binom{k+2}{k}\sum_{i=0}^k (-1)^i\binom{k}{i}\binom{n-\frac i2-2}{n-j}
\]
where $\binom{a}{b}=a(a-1)\dots(a-b+1)/b!$.
\end{pro} 
This model of urns has also been studied by Janson \cite{SJ}; Theorem 1.3 in \cite{SJ} gives the asymptotic behavior of urns under these dynamics, depending on the initial conditions. The discussion given in Section 3.1 of \cite{SJ} shows that the asymptotic behavior of $D_j(n)$ is quite difficult to describe. 
One may use \eref{evold} to see that $`E(D_j^{n+1} ~|~D_j^n)=D_j^n(1+\frac{1}{2n-1})$ to show that $(M_j^{n})_{n\geq j}$ defined by 
\[M_j^{n}=D_j^n / u_n\]
is a ${\cal F}_n$ martingale, where ${\cal F}_n=\sigma(D_j^{k}, j\leq k\leq n )$
for any sequence $u_n$ such that $u_{n+1}=u_n (2n)/(2n-1)$. This allows to see that 
\[`E(D^n_j)= `E(D^j_j)\prod_{k=j}^{n-1} (2k)/(2k-1)=3 \prod_{k=j}^{n-1} (2k)/(2k-1).\]
This indicates that for a fixed $j$ the expectation $`E(D^n_j)$ grows as $\sqrt{n}$. 
Some other regimes may be obtained: for $t\in(0,1)$, $`E(D^n_{nt})\to 3(1-t)^{-1}$ when $n$ goes to $+\infty$. 
(We recall that any triangulation with $2n$ faces has $3n$ edges and $n+2$ nodes; hence the mean degree of a node in $6n/(n+2)$ in any triangulation).

\section{Two families of increasing quadrangulations}
\label{allquad}
We present here two families of quadrangulations. The first one, quite natural, resists to our investigations. The second one, that may appear to be quite unnatural, is in fact very analogous to stack-triangulations, and is studied with the same tools.

\subsection{A first model of growing quadrangulations}
\label{qua-dur}
This is the simplest model, and we present it rapidly: starting from a rooted square, choose a finite face $f=ABCD$ and a diagonal $AC$ or $BD$. Then add inside $f$, a node $x$ and the two edges $Ax$ and $xC$ or the two edges $Bx$ and $xD$. The set $\4{k}'$ is then the set of quadrangulations with $n$ bounded faces reached by this procedure starting with the rooted square (formally define a growing procedure $\Phi_4$, similar to $\Phi$ of Section \ref{def}, using $\4{k}'{}^{\bullet}=\{(m,f,\alpha) ~|~ m\in \4{k}',  f\in F^\circ(m), \alpha\in \{0,1\}\}$ the rooted quadrangulations from $\4{k}'$ with a distinguished finite face marked with 0 or 1, and add in $f$ a pair of edges or the other one according to $\alpha$). 
\medskip

There is again some bijections between $\4{k}'$ and some set of trees, but we were unable to define on the corresponding trees a device allowing to study the distance in the maps (under the uniform distribution, as well as under the distribution induced by the construction when both $f$ and $\alpha$ are iteratively uniformly chosen).  We conjecture that they behave asymptotically in terms of metric spaces as triangulations under $\qt_{2k}$ and $\ut_{2k}$ up to some normalizing constant. \par

We describe below a bijection between $\4{k}'$ and the set of trees having no nodes having only one child. There exists also a bijection with Schr\"oder trees (trees where the nodes have 0,1 or 2 children) with $k$ internal nodes. 
\begin{pro}For any $k\geq 2$, there exists a bijection $\Psi_k$ between $\4{k}'$ and the set of trees having $k$ leaves, no nodes of outdegree 1 and with a root marked 0 or 1.
\end{pro}
 For $k=1$, $\4{k}'=\{s\}$ the rooted square and in this case we may set $\Psi(s)=\{\varnothing\}$, the tree reduced to a (non marked) leaf.\\
\proof
Assume that $k\geq 2$. Split $\4{k}'$ into two subsets $\4{k,0}'$ and $\4{k,1}'$, letting $\4{k,0}'$ contains the maps $m$ with exterior faces $ABCD$ rooted in $AB$ containing an internal node $x$ and the two edges $Ax$ and $xC$, and $\4{k,1}'$ those containing an internal node $x$ and the two edges $Bx$ and $xD$ (notice that $m$ cannot contain at the same time an internal node $x$ and $Ax$ and $xB$). It is easy to see that the rotation of $\pi/2$ is a bijection between $\4{k,0}'$ and $\4{k,1}'$. We then focus on $\4{k,0}'$ and explain the bijection between $\4{k,0}'$ and the set of trees having no nodes of outdegree 1 and $k$ leaves.
Let $x_1,\dots,x_j$ be the $j\geq 1$ internal points of $m$, adjacent to $A$ and $C$. These points (if properly labeled) define $j+1$ submaps $m_1,\dots,m_{i+j}$ of $m$ with border $Ax_iCx_{i+1}$ for $i=0$ to $j$ where $B=x_0$ and $D=x_{j+1}$. We then build $t=\Psi_k(m)$ by sending $m$ onto the root of $t$, 
and $m_i$ to the $i$th child of $m$.
Each of the submaps $m_i$ can also be decomposed in the same way except that by maximality of the set $\{x_1,\dots,x_j\}$, the face $Ax_iBx_{i+1}$ is either empty or contains an internal node $y$ adjacent to $x_i$ and $x_{i+1}$. The coloring of the nodes (except) the root is then useless.  ~$\Box$\medskip

\subsection{A family of stack-quadrangulations}
\label{fsq}
The construction presented here is very similar to the construction of
stack-triangulations; some details will be skipped when the analogy with
them will be clear enough. The difference with the model of quadrangulations of Section \ref{qua-dur} is that given a face $f=ABCD$, only a suitable choice of pair of edges (either $(Ax,xC)$ or $(Bx,xD)$) will be allowed.

This choice amounts to forbidding double ``parallel'' pair of edges of the type $(Ax,xC)$ and $(Ax',x'C)$.
\par

Formally,  set first $\4{1}=\{s\}$ where $s$ is the unique rooted square. There is also a unique element in  $\4{2}$ obtained as follows. Label by $ABCD$ the vertices of $s$, such that $(A,B)$ is the root of $s$. To get the unique element of $\4{2}$, draw $s$ in the plane, add in the bounded face of $s$ a node $x$ and then the two edges $(Ax)$ and $(xC)$ in this face. \par
\begin{figure}[htbp]
\psfrag{A}{A}
\psfrag{B}{B}
\psfrag{C}{C}
\psfrag{D}{D}
\psfrag{x}{x}
\centerline{\includegraphics[height= 2.5 cm]{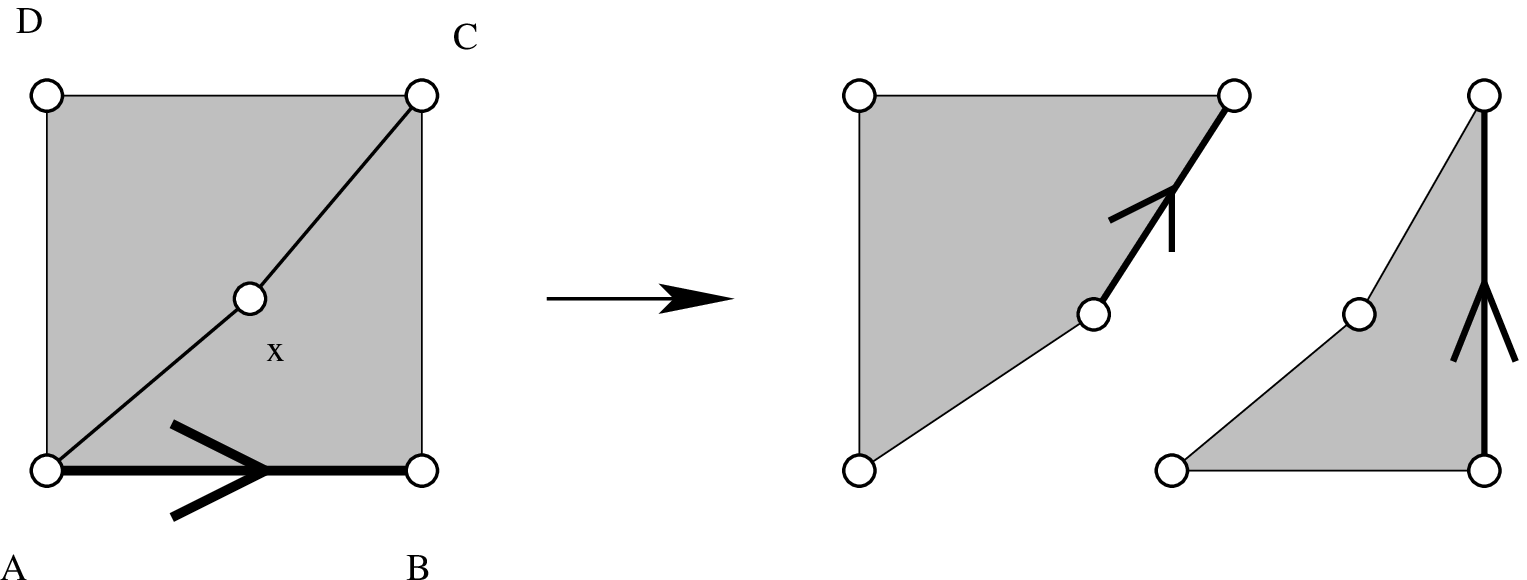}}
\captionn{\label{decompote} The decomposition is well defined thanks to the uniqueness of a node $x$ adjacent to both $A$ and $C$.}
\end{figure}

 
We define now $\4{k}$ recursively asking to the maps $m$ with border $ABCD$ and rooted in $(A,B)$ to have the following properties. If $k\geq 1$ there exists a unique  node $x$ in the map $m$, such that $Ax$ and $xC$ are edges of $m$. Moreover the submaps $m_1$ and $m_2$ of $m$ with respective borders $AxCD$ (rerooted in $(x,C)$) and $ABCx$ (rerooted in $(B,C))$ belong both to the sets $\cup_{j< k}\4{j}$, more precisely $(m_1,m_2)\in \cup_{j=1}^{k-1} \4{j}\times\4{k-j}$ (see an illustration on Figure \ref{decompote}).\par

This rerooting operation corresponds to distinguish a diagonal in each face (once for all) on which the following insertion inside this face, if any, will take place. \par 
\par

Any maps belonging to $\4{k}$ is a rooted quadrangulation having $k$ internal faces. There exists again a canonical drawing of these maps, where the border $ABCD$ (rooted in $(A,B)$) of the quadrangulations is sent on a fixed square of the plane, and where, when it exists, the unique node $x$ adjacent to both $A$ and $C$ is sent of the center of mass of $ABCD$, the construction being continued recursively in the submaps $m_1$ and $m_2$ (the edges are straight lines).\par

There exists also a sequential construction of this model, more suitable to define the distribution of interest. 

This is very similar to the case of triangulations treated in Propositions \ref{yop} and \ref{yopp}.

\subsubsection{Sequential construction of $\4{k}$}
\label{scq}
We introduce a labeling of the nodes of $\4{k}$ by some integers. The idea here is double. 
This labeling will distinguish the right diagonal where we will insert pair of edges, and also, will be used to count the number of histories leading to a given map. A labeled map may be viewed as a pair $(m,l)$ where $m$ is an unlabeled map and $l$ an application from $V(m)$ onto the set of integers.\par
We then consider $\4{k}^\ell$ be the set of quadrangulation having $k$ internal faces and where the vertices are labeled as follows. 
First $\4{1}^\ell$ contains the unique labeled rooted map $(s,l)$ with vertices $ABCD$ (rooted in $(A,B)$) and labeled by
\[l(A)=4, l(B)=3, l(C)=2, l(D)=1.\]
Assume now that $\4{k}^\ell$ is a set of quadrangulations with $k$ internal faces  (and thus $k+3$ vertices), where the vertices are labeled by different integers from $\{1,\dots,k+3\}$.  
To construct $\4{k+1}^\ell$ from $\4{k}^l$ we consider an application 
$\Phi^\ell_4$ from $\4{k}^{\ell,\bullet}=\{((m,l),f) ~|~ m\in\4{k}^l, f\in m^\circ\}$ such that:
to obtain $\Phi^l_4((m,l),f)$, draw the label map $m$ in the plane; denote by $ABCD$ the vertices of $f$, such that $A$ has the largest label (and thus $C$ is at the opposite diagonal of $A$ in $f$).  Add a point $x$ labeled $k+4$ in $f$ and the two edges $Ax$ and $xC$ in $f$. The obtained labeled map is $\Phi^l_4((m,l),f)$.

We denote by $\4{k+1}^\ell$ the set $\Phi^\ell_4(\4{k}^{\ell,\bullet})$. \par

We call $\pi_k$ (or more simply $\pi$) the function
\[\app{\pi_k}{\4{k}^\ell}{\4{k}}{(m,l)}{m}\]
the canonical surjection from $\4{k}^\ell$ onto $\4{k}$; this is simply the application that erases the labels of a labeled map. This definition hides a property, since the set $\4{k}$ has been defined in the beginning of Section \ref{fsq}. The proof of the equality of the sets $\pi(\4{k}^\ell)$ and $\4{k}$ is a consequence of the binary decomposition of both object according to the distinguished diagonal.

\begin{figure}[htbp]\psfrag{x}{$x$}
\psfrag{A}{$A$}\psfrag{1}{1}\psfrag{2}{2}\psfrag{3}{3}\psfrag{4}{4}
\psfrag{B}{$B$}\psfrag{5}{5}\psfrag{6}{6}\psfrag{7}{7}\psfrag{8}{8}
\psfrag{C}{$C$}
\psfrag{D}{$D$}
\centerline{\includegraphics[height= 2.6 cm]{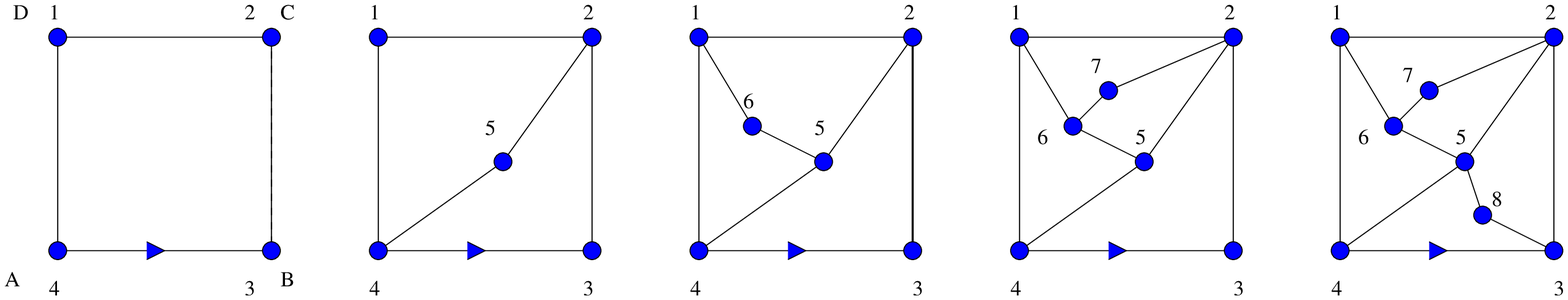}}
\captionn{\label{fig34}A sequence of quadrangulations obtained by successive insertions of pair of edges. }
\end{figure}

Consider a labeled map $(m_k,l_k)\in\4{k}^{\ell}$ for some $k\geq 2$. There exists a unique map $(m_{k-1},l_{k-1})\in \4{k-1}^{\ell}$ such that $(m_k,l_k)=\Phi^\ell_4(m_{k-1},l_{k-1})$. It is obtained from $(m_k,l_k)$ by the suppression of the node with largest label together with the two edges that are incident to this node. Hence, each map 
$(m_k,l_k)$ characterizes uniquely a legal history of $m_k=\pi(m_k,l_k)$. We mean by ``legal'' that for any $j$, $m_{i+1}$ is obtained from $m_i$ by the insertion of two edges, and for any $i$, $m_i$ is in $\4{i}$. From now on, we will make a misuse of language and confound the histories of a stack-quadrangulation $m_k\in \4{k}$ and $\pi^{-1}(m_k)$.

We denote by $\uq_{k}$ the uniform distribution on $\4{k}$
and as we did for triangulations in Section \ref{tw}. 
We denote by $\qq_{k}$ the distribution of $\pi(M_k,l_k)$ when $M_{i+1}=
\Phi_4^\ell((M_i,l_i),F_i)$, where $M_1$ is the only element of $\4{1}$ and where $F_i$ is chosen
uniformly among the internal faces of $M_i$ (all the $F_i$ are independent).
The support of $\qq_{k}$ is the set $\4{k}$, and one may check that $\qq_{k}\neq\uq_{k}$ for $k\geq 4$.

\subsection*{The function $\Gamma'$}

As in Section \ref{fond-bij}, we define a function $\Gamma'$ to express the
distance between any pair of nodes $u$ and $v$ in a stack-quadrangulation
$m$ in terms of a tree associated bijectively to this map.
Let 
\[
W_{1,2} = \{12,21\}^{\star}\cdot\{11,22\}\cdot\{1,2\},
\]
be the set of words on $\Sigma_2=\{1,2\}$, beginning with any number of occurrences of
$12$ or $21$, followed by $11$ or $22$, then by a $1$ or a $2$. Notice  that all the words of $W_{1,2}$ have a odd length.
For example $u=12\, 21\, 21\,  11\, 2 \in W_{1,2}$.

Let $u=u_{1}\ldots u_k$ be a word on the alphabet $\Sigma_2$. Define
$\tau_1(u):=0$ and for $j\geq 2$,
\[
\tau_j(u) :=\inf\{i ~|~ i \ge \tau_{j-1}(u) \text{ such that }
u_{1+\tau_{j-1}(u)}\dots u_i\in W_{1,2} \}.
\]
This amounts to decomposing $u$ into subwords belonging to $W_{1,2}$.
We denote by $\tilde{\Gamma'}(u)=\max\{i ~|~ \tau_i(u)\leq |u|\}$, then
$u=u_1\ldots u_{\tau_{\tilde{\Gamma'}(u)}(u)}\tilde{u}$, where $\tilde{u}\notin W_{1,2}$.
Lastly we define $\Gamma'(u)$ as
\[
\Gamma'(u)=\tilde{\Gamma'}(u) +
\begin{cases}
  0 & \text{if }|\tilde{u}|\text{ is even and } \tilde{u} \text{ does not end with }11
  \text{ or }22\\
  1 &\text{otherwise}
\end{cases}
\]
Further, for two words $u=w a_1 \ldots a_k$ and $v=w b_1 \ldots b_l$ (with
$a_1\neq b_1$), set as in the triangulation case $\Gamma'(u,v)=\Gamma'(a_1\ldots a_k) + \Gamma' (b_1\ldots b_l)$.\par
We now give a proposition similar to Proposition \ref{yop} for stack-quadrangulations.
\begin{pro}\label{prop:bijquad}
For any $K\geq 1$, there exists a bijection
\[
\app{\Psi_K^{`4}}{\4{K}}{\Tbin_{2K-1}}{m}{t:=\Psi_K^{`4}(m)}\]
such that~:\\
$(i)$ $(a)$ each internal node $u$ of $m$ corresponds bijectively to an internal node $u'$ of $t$.\\
$(b)$ Each leaf of $t$ corresponds bijectively to a finite quadrangular face of $m$. \\
$(ii)$ For any $u$ internal node of $m$, $\Gamma'(u')=d_m(root,u).$ \\
$(ii')$ For any $u$ and $v$ internal nodes of $m$
\begin{equation}\label{eq:degquad}
\l|d_{m}(u,v)-\Gamma'(u',v')\r|\leq 4.
\end{equation}
\noindent$(iii)$ Let $u$ be an internal node of $m$. We have
\[\deg_m(u)=2+\#\{v'\in t^{\circ} ~|~ v'=u'w',|w'|\geq 2, w'\in \{12,21\}^{\star}\},\]
\end{pro}
The existence of a bijection between $\4{K}$ and $\Tbin_{2K-1}$ comes from
the recursive decomposition of a stack-quadrangulation along the first 
pair of edges inserted (which can be determined at any time since there  is  a unique node $x$ adjacent to both $A$ and $C$ in any $m \in \4{K}$, for $K\geq 2$). 

\begin{figure}[htbp]
\psfrag{A}{$A$}
\psfrag{B}{$B$}
\psfrag{C}{$C$}
\psfrag{D}{$D$}
\centerline{\includegraphics[height= 5.5 cm]{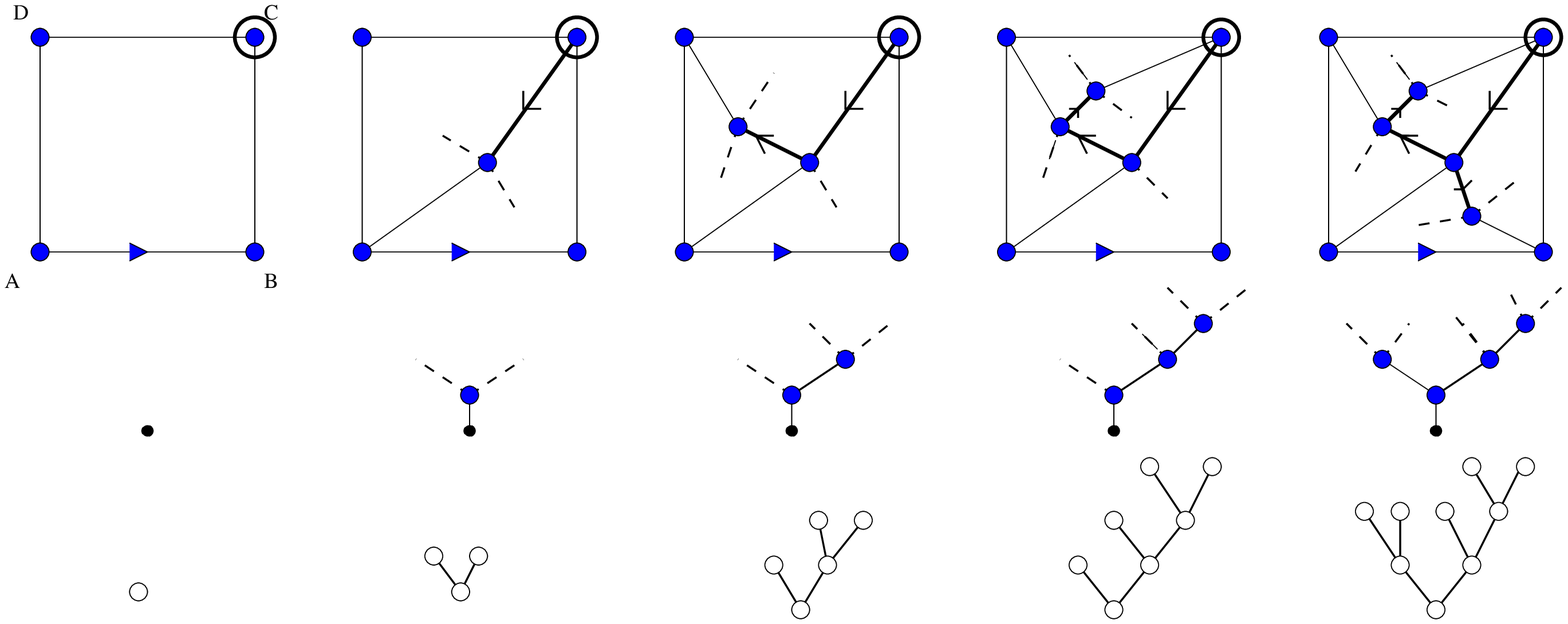}}
\captionn{\label{fig:dec_quad}A sequence of quadrangulations obtained by successive insertions of pair of edges. }
\end{figure}

\proof  The proof of this Proposition is very similar to that of Proposition \ref{yopp}.
We only sketch the main lines. First, the maps in $\4{K}$ own also a canonical drawing as said above. 
We propose a bijection that does not follow the decomposition provided in Figure \ref{decompote}, but which is illustrated in Figure \ref{fig:dec_quad}. 

Hence, we start from the square $(A,B,C,D)$ rooted in $(A,B)$. We associate with any stack-quadrangulation a binary tree as represented on Figure \ref{fig:dec_quad}. Again, the possibility to think in terms of canonical maps and faces, allow to see the consistence and robustness of the sequential approach represented on the illustration. If $u$ is associated with a face, then  $u1$ (resp. $u2$) corresponds to the face situated on the left (resp. on the right) of this oriented edge (see Figure \ref{fig:dec_quad}).

\par
To prove properties $(ii)$ and $(iii)$ we introduce a notion of type of faces in a stack triangulation (or type of a node in the corresponding tree) as in the proof of Proposition \ref{yopp}.
 For any face $f=(A,B,C,D)$ in $m$ such that
$O(f)=(A,B)$, we set:
\[
\text{type}(A,B,C,D):=(d_m(E_0,A),d_m(E_0,B),d_m(E_0,C),d_m(E_0,D))
\]
the 4-tuple of the distance of $A$, $B$, $C$ and $D$ to the root vertex of
$m$. 
It is well known that in a quadrangulation, the type of any face is $(i,i+1,i,i+1)$ or $(i,i+1,i+2,i+1)$, for some $i$, or a circular permutation of this. 

As the types of the faces arising in the construction are not modified by the insertions of new edges, we mark any node of $t=\Psi^{`4}(m)$ with the type of the corresponding face. It is then easy to check that for $u'$ an internal node of $t$ with type$(u')=(a,b,c,d)$, we have $d_m(u,E_0)=1+b\wedge d$ and 
\begin{equation}\label{eq:type_qua}
\left\{
\begin{array}{cccccl}
\type(u'1)=(&b,&1+b\wedge d,&d,&a&),\\
\type(u'2)=(&b,&1+b\wedge d,&d,&c&),\\
\end{array}\right.
\end{equation}
 Property $(ii)$ follows directly from (\ref{eq:type_qua}) using the fact that
type$(\varnothing)=(1,2,1,0)$. Properties $(ii')$ and $(iii)$ are deduced directly
by the same arguments as for triangulations.\hfill $\Box$

\subsubsection{Asymptotic behavior of the quadrangulations}
First, we state a Lemma analogous to Lemma \ref{cou}, that can be proved similarly except that here  $\tau_{i}-\tau_{i-1}\sim 1+\Geo(1/2)$ and thus has mean 3 (we also use that for any $u\in \{1,2\} ^{\star}$, $|\Gamma'(u)-\tilde{\Gamma'}(u)|\leq 1$).
\begin{lem}\label{lem:asquad}
Let $(X_i)_{i\geq 1}$ be a sequence of i.i.d. random variables taking their values in $\Sigma_2=\{1,2\}$ and let $W_n$ be the word $(X_1,\dots,X_n)$.\\
$(i)$ $n^{-1}{\Gamma'(W_n)}\as\Gamma'_{`4}$
where
\begin{equation}
\Gamma'_{`4}:=1/3
\end{equation}
$(ii)$ $`P( |\Gamma'(W_n)-n\Gamma'_{`4}|\geq n^{1/2+u})\to 0$ for any $u>0$.
\end{lem}

We are now in position to state the main theorem of this part. We need to examine first the weak limit of binary trees. Denote by $P^{\bin}_{2n+1}$ be the uniform distribution on the set of binary trees with $2n+1$ nodes. 
This time $P^{\bin}_{\infty}$ is the distribution of a random infinite tree, build around an infinite line of descent $L^\bin_{\infty}=(X(j), j\geq 0)$, where $(X(j))$ is a sequence of i.i.d. r.v. uniformly distributed on $\Sigma_2=\{1,2\}$ on the neighbors of which are grafted critical GW trees with offspring distribution $\nu_\bin=(1/2)(\delta_0+\delta_2)$.
We sum up in the following Proposition, the results concerning the convergence of trees under  $P^{\bin}_{2n+1}$.
\begin{pro}\label{crtbin}
$(i)$ The following convergence holds for the GH topology. Under $P^{\bin}_{2n+1}$,
\[\l(T, \frac{d_T}{\sqrt{2n}}\r)\dd ({\cal T}_{2\se},d_{2\se}).\]
$(ii)$ When $n\to+\infty$, $P^{\bin}_{2n+1}$ converges weakly to $P^{\bin}_\infty$ for the topology of local convergence.
\end{pro}
The first point is due to Aldous \cite{ALD}, and the second point (very similar to Proposition \ref{loctree}) is also due to Gillet \cite{FG}.

Again, all the results of Section \ref{rp} may be extended in the binary case, as well as the construction of the infinite map $(m_\infty^\square)$ in a way similar to $m_\infty$, the  limit of triangulations for the local convergence. We can then prove, following the lines of the triangulation case
\begin{theo}\label{loctri2} $(i)$ Under $\uq_{n}$, $(m_n)$ converges in distribution to $m_{\infty}^\square$ for the topology of local convergence.\\
\label{youp2}
$(ii)$ Under $\uq_{n}$,   
\[\l(m_n,\frac{D_{m_n}}{\Gamma'_\square\sqrt{2n}}\r)\dd ({\cal T}_{2\se}, d_{2\se}),\]
for the Gromov-Hausdorff topology on compact metric spaces.
\end{theo}

Now, the asymptotic behavior of maps under $\qq_k$ are studied again thanks to trees under $\qbin_{2K-1}:=\qq_{K}\circ (\Psi_K^{`4})^{-1}$ the corresponding distribution on trees. This distribution on $\Tbin_{2K-1}$ is famous in the literature since it corresponds to the distribution of binary search trees. Indeed the insertion in the map $m$, corresponds to an uniform choice of a leaf in the tree $\Psi^{`4}(m)$ and its transformation into an internal node having two children. Again, using the same tools as those used to treat the asymptotic behavior of trees under $\qter$ (in particular, here the fragmentation is binary, and $Y^u\sur{=}{d}(U,1-U)$ where $U$ is uniform in $[0,1]$), we get the following proposition.
\begin{pro}\label{dic2}
Let $\bt$ a random tree under the distribution $\qbin_{2K+1}$. Let ${\bf u}$ and ${\bf v}$ be two internal nodes chosen equally likely and independently among the internal nodes of $\bt$, and let ${\bf w}={\bf u}\wedge {\bf v}$.\\
1) We have
\[\l(4\log n\r)^{-1/2}\l(|{\bf u}|-4\log n,|{\bf v}|-4\log n\r)\dd (N_1,N_2)\]
where $N_1$ and $N_2$ are independent centered Gaussian r.v. with variance 1.\\
2) Let ${\bf a}, {\bf b}\in \{1,2\}$, with ${\bf a}\neq {\bf b}$ and  ${\bf u}^\star, {\bf v}^\star$ the (unique) words such that
\[{\bf u}={\bf wau}^\star \textrm{ and } {\bf v}={\bf wbv}^\star.\]
Conditionally to $(|{\bf u}^\star|,|{\bf v}^\star|)$ (their lengths) ${\bf u}^\star$ and ${\bf v}^\star$ are independent random words composed with $|{\bf u}^\star|$ and $|{\bf v}^\star|$ independent letters uniformly distributed in $\Sigma_2=\{1,2\}$.
\end{pro}
The interested reader may find in Mahmoud \& Neininger \cite[Theorem 2]{MN} a different proof of the first assertion, the second one, once again being a consequence of the symmetries of this class of random trees.

Similarly to Theorem \ref{metconv}, we obtain
\begin{theo}\label{metconvquad} Let $M_n$ be a stack-quadrangulation under $\qq_{2n}$. Let $k\in \mathbb{N}$ and ${\bf v}_1,\dots,{\bf v}_k$ be $k$ nodes of $M_n$ chosen independently and uniformly among the internal nodes of $M_n$. We have
\[\l(\frac{D_{M_n}({\bf v}_i,{\bf v}_j)}{4\Gamma'_\square\log n}\r)_{(i,j)\in\{1,\dots,k\}^2}\proba \l(1_{i\neq j}\r)_{(i,j)\in\{1,\dots,k\}^2}\] the matrix of the discrete distance on a set of $k$ points.
\end{theo}

\section{Appendix}
\label{ap}
\subsection{Proof of the Theorems of Section \ref{asGH}}
\label{rp}
The  aim of this section is to prove Theorem \ref{paramversion}. 
Our study of the distance in a stack-triangulation $m_n$ passes via the study of the function $\Gamma$ on the tree $T=\Psi_n^{`3}(m_n)$. 
Let $w(r)$ be the $r$th internal node of $T$ according to the LO ($w(0)$ is the root), and $u(r)$ be the $r$th internal node of $m$ (the image of $w(r)$ as explained in Proposition \ref{yop}). For any $r$ and $s$, 
\begin{equation}\label{mou}
\l|d_{m}(u(r),u(s))-\Gamma(w(r),w(s))\r|\leq 4.
\end{equation}

\begin{lem}Under $\ut_{2n}$,  the family $\l(\l(\frac{d_{m_n}(ns,nt)}{\Gamma_\triangle\sqrt{3n/2}}\r)_{(s,t)\in[0,1]^2}\r)_n$ is tight on $C[0,1]^2$.
\end{lem}

\proof 
We claim first that under $\uter_{3n+1}$, the family  $\l(n^{-1/2}d_{T^\circ}(ns,nt)\r)_n$ is tight in $C[0,1]^2$, where $d_{T^\circ}(k,j)=d_{T^\circ}(w(k),w(j))$, is the (reparametrization of the) restriction of the distance in $T$ 
on its set of internal nodes, and where  $d_{T^\circ}$ is smoothly interpolated as explained below Theorem \ref{youp}.
Indeed, let $(H^\circ(k))_{k=0,\cdots,n-1}$ where $H^\circ(k)=|w(k)|$ be the height process of the internal nodes of $T$ (interpolated between integer points). According to Marckert \& Mokkadem \cite[Corollary 5]{MMexc} (and using that the height process of ternary trees coincides with the height process restricted to node with outdegree 3), 
\begin{equation}\label{convhp}
\l(\frac{H^\circ({nt})}{\sqrt{3n/2}}\r)_{t\in[0,1]}\dd \l({2}\se_t\r)_{t\in[0,1]}.
\end{equation}
Using that for $i\leq j$, 
\[\l|d_T(w(i),w(j))-(H^\circ(i)+H^\circ(j)-2 \min_{k\in \cro{i,j}} {H^\circ}(k))\r|\leq 2\]
we get that 
\[\l(\frac{d_{T^\circ}(ns,nt)}{\sqrt{3n/2}}\r)_{s,t\in[0,1]}\dd \l(d_{2\se}(s,t)\r)_{s,t\in[0,1]}.\]
where the convergence holds in $C[0,1]^2$. This is just a consequence of the continuity of the application $f\mapsto \l[(s,t)\to f(s)+f(t)-2\min_{u\in[s,t]}f(u)\r]$ from $C[0,1]$ onto $C[0,1]^2$.
We deduce from this that the sequence 
$\l(\l(\frac{d_{T^\circ}(ns,nt)}{\sqrt{3n/2}}\r)_{s,t\in[0,1]}\r)_n$ is tight and
by \eref{mou} and the trivial bound 
$\Gamma(u,v)\leq d_{T^\circ}(u,v)$ for any $u$ and $v\in T^\circ,$
\[\frac{d_{m_n}(ns,nt)}{\sqrt{3n/2}}\leq \frac{d_{T^\circ}(ns,nt)}{\sqrt{3n/2}}+4n^{-1/2}\]
and thus the Lemma holds true.~$\Box$.

The convergence of the finite dimensional distributions in Theorem \ref{paramversion} is a consequence of the following stronger result.
\begin{pro}\label{fond}
Let $0\leq s<t\leq 1$. 
When $n$ goes to $+\infty$, under $\ut_{3n+1}$
\[\l|\frac{d_{m_n}(\floor{ns},\floor{nt})}{\Gamma_\triangle\sqrt{3n/2}}-\frac{d_{T\circ}(\floor{ns},\floor{nt})}{\sqrt{3n/2}}\r|\proba 0.\]
\end{pro}

To prove this Proposition we need to control precisely $\Gamma(w(ns),w(nt))$; we will show that this quantity is at the first order, and with a probability close to 1, equal to $\Gamma_\triangle d_{T^\circ}(ns,nt)$. This part is largely inspired by the methods developed in a work of the second author \cite{MJF}. \par
 
We focus only on the case $s,t$ fixed in $(0,1)$ and $s<t$ (which is the most difficult case). In the following we write $ns$ and $nt$ instead of $\floor{ns}$ and $\floor{nt}$.
 Consider  $\cw_{ns,nt}=w(ns)\wedge w(nt)$, and write
\begin{equation}
\label{dec-word}
w({{ns}})=\cw_{ns,nt} l_0 l_{ns,nt}~~ \textrm{ and }~~w({nt})= \cw_{ns,nt} r_0 r_{ns,nt},
\end{equation}
where $l_0\neq r_0$  (the letters $l$ and $r$ refer to ``left'' and ``right'').

For compactness of notation, set
\be
\Dec(n)&:=&(W_1,W_2,W_3,H_1,H_2,H_3,L,R)\\
&:=&(\cw_{ns,nt},l_{ns,nt},r_{ns,nt},
|\cw_{ns,nt}|,|l_{ns,nt}|,|r_{ns,nt}|,l_0,r_0),
\ee
$\Dec$ standing for "decomposition". Even if not recalled in the statements,
these variables are considered as random variables under $`P_{3n+1}^{\ter}$. Let now $\widetilde{\Dec}$ be the random variable defined by 
\[\widetilde{\Dec}(n):=
(\tilde{W_1},\tilde{W_2},\tilde{W_3},H_1,H_2,H_3,\tilde{L},\tilde{R})\]
such that, 
conditionally on $(H_1,H_2,H_3)=(h_1,h_2,h_3)$, the random variables $\tilde{W_1},\tilde{W_2},\tilde{W_3},\tilde{L},\tilde{R}$ are independent and defined by:\\
-- for each $i\in\{1,2,3\}$, $ \tilde{W_i}$ is a word with $h_i$ i.i.d. letters, uniformly chosen in $\{1,2,3\}$,\\
-- the variable $(\tilde{L},\tilde{R})$ is a random variable uniform in $I_3=\{(1,2),(1,3),(2,3)\}$.

\begin{defi}
Let $(Y_1,Y_2,\dots)$ and $(X_1,X_2,\dots)$ be two sequences of r.v. taking their values in a Polish space $S$. We say that $`P_{X_n}/`P_{Y_n}\sur{\to}{\star} 1$ or  $X_n{\di}_{\!\star}~ Y_n\to 1$ if for any $`e>0$ there exists a measurable set $A_n^{`e}$ and a measurable function $f_n^{`e}: A_n^{`e}\mapsto `R$ satisfying $`P_{X_n}=f_n^{`e} `P_{Y_n}$ on $A_n^{`e}$, such that $\sup_{x\in A_n^{`e}}|f_n^{`e}(x)-1|\sous{\to}{n}0$ and such that $`P_{Y_n}(A_n^{`e})\geq 1-`e$  for $n$ large enough.
\end{defi}

The main step in the proof of Proposition \ref{fond} is the following Proposition. 

\begin{pro}\label{cvradon} When $n\to +\infty$, $\Dec(n){\di}_{\!\star} \widetilde{Dec}(n)\to 1.$
\end{pro}

Assume that this proposition holds true and let us end the proof of Proposition \ref{fond}. 
The  following lemma (proved in \cite[Lemma 16]{MJF}\footnote{In\cite[Lemma 16]{MJF} the function $g_n$ is assumed to be continuous, but only the measurability is needed}) allows to compare the limiting behavior of $\Dec(n)$ and $\widetilde{\Dec}(n)$. 
\begin{lem}\label{dr3}Assume that ${X_n}{\di}_{\!\star} Y_n\to 1$ then:\\
$\bullet$ If $Y_n\dd Y$ then $X_n\dd Y$.\\
$\bullet$ Let $(g_n)$ be a sequence of measurable functions from $S$ into a Polish space $S'$. If ${X_n}{\di}_{\!\star}~ Y_n\to 1$ then $g_n(X_n){\di}_{\!\star} ~g_n(Y_n)\to 1$
\end{lem}
\prooff{Proposition \ref{fond}} From  Proposition \ref{cvradon} and Lemma \ref{dr3}, we deduce
\[(H_2,H_3,W_2,W_3){\di}_{\!\star}(H_2,H_3,\tilde{W}_2,\tilde{W}_3)\to 1.\] 
Since $(3n/2)^{-1/2}\l({H_2,H_3,\Gamma(\tilde W_2),\Gamma(\tilde W_3)}\r)$ converges in distribution to 
\begin{equation}\label{limit}
\l(2\se_s-m_{2\se}(s,t),2\se_t-m_{2\se}(s,t),\Gamma_\triangle(2\se_s-m_{2\se}(s,t)),\Gamma_\triangle(2\se_t-m_{2\se}(s,t))\r)
\end{equation}
thanks to Lemmas \ref{convhp}, \ref{cou}  (and also Lemma \ref{qd} below which ensures that $H_i \in [M^{-1},M]\sqrt{n}$ with probability arbitrary close to 1, if $M$ is chosen large enough, leading to a legal using of Lemma \ref{cou}). 
We then deduce by the first assertion of Lemma \ref{dr3} that
$(3n/2)^{-1/2}\l(H_2,H_3,\Gamma(W_2),\Gamma(W_3)\r)$ converges also in distribution to the random variable described in \eref{limit}.
In particular this implies
\[n^{-1/2}\l|\Gamma_\triangle d_{T}(w(ns),w(nt))-\Gamma(w(ns),w(ns))\r|\proba 0.~~~\Box\]

It only remains to show Proposition \ref{cvradon}. The absolute continuity $`P_{\Dec(n)}\prec `P_{\widetilde{\Dec}(n)}$ comes from the inclusion of the (discrete) support of $\Dec(n)$ in that of $\widetilde{\Dec}(n)$. \medskip

For any word $w=w_1\dots w_k$ with letters in $\{1,2,3\}$ define
\[N_1(w)=\sum_{j=1}^k (w_i-1) \textrm{ and }N_2(w)=\sum_{j=1}^k (3-w_i). \]
Seeing $w$ as a node in a tree, $N_1(w)$ and $N_2(w)$ give the number of nodes at distance 1 on the left (resp. on the right) of the branch $\cro{\varnothing, w}$. 
Set 
\be
A_{n,M}&=&\{(w_1,w_2,w_3,h_1, h_2, h_3,l,r) ~|~ h_1,h_2,h_3 \in\sqrt{n}[M^{-1},M],\\
&& (w_1,w_2,w_3)\in J_{h_1}\times J_{h_2} \times J_{h_3}, (l,r)\in I_3\},
\ee
where for any $h>0$,
$J_{h}=\l\{\sa \in\Sigma_3^h ~|~ (N_1(\sa),N_2(\sa))\in \Big[h- h^{2/3},h+ h^{2/3}\Big]^2\r\}$.
\begin{lem}\label{qd} For any $`e>0$, there exists $M>0$ such that for $n$ large enough
\[`P_n(\widetilde{Dec(n)}\in A_{n,M})\geq 1-`e.\]
\end{lem}
\proof 
First, by the Skohorod's representation theorem (see e.g.\ \cite[Theorem 4.30]{KAL}) there exists a probability space where the convergence of the rescaled height process to $2\se$ (as stated in \eref{convhp}) is an a.s. convergence. On this space $(3n/2)^{-1/2}(H_1,H_2,H_3)$ converges a.s. to $(m_{2\se}(s,t),2\se_s,2\se_t)$. Since  a.s. $m_{2\se}(s,t)<\min(e_s,e_t)$,  thus
\begin{equation}\label{dist}
\liminf `P^\ter_{3n+1}(H_i\in [M^{-1},M]\sqrt{n}, i\in\{1,2,3\})\geq 1-`e \textrm{~~~for }M \textrm{ large enough}.
\end{equation}
Let $W[h]$ be a random word with $h$ i.i.d. letters uniform in $\Sigma_3$. 
 For  $h\in\mathbb{N}$, by symmetry $N_1(W[h])$ and $N_2(W[h])$ have the same law, and there exists $c_1>0,c_2>0$, s.t
\[`P\l(W[h]\notin J_{h}\r)\leq c_1\exp(-c_2\,h^{1/3}).\]
Indeed  the number $x_i$ of letters $i$ in $W[h]$ is binomial $B(h,1/3)$ distributed, and the Hoeffding inequality leads easily to this result ($N_1(h)= x_2+2x_3$ which is in mean $h/3+2h/3=h$). ~$\hfill\Box$\medskip

\begin{figure}[htbp]
\psfrag{v}{$\varnothing$}
\psfrag{w}{$\check{u}_{ns,nt}$}
\psfrag{s}{${u}({ns})$}
\psfrag{t}{${u}({nt})$}
\psfrag{F_1}{$F_1$}
\psfrag{F_2}{$F_2$}
\psfrag{F_3}{$F_3$}
\centerline{\includegraphics[height= 4.4 cm]{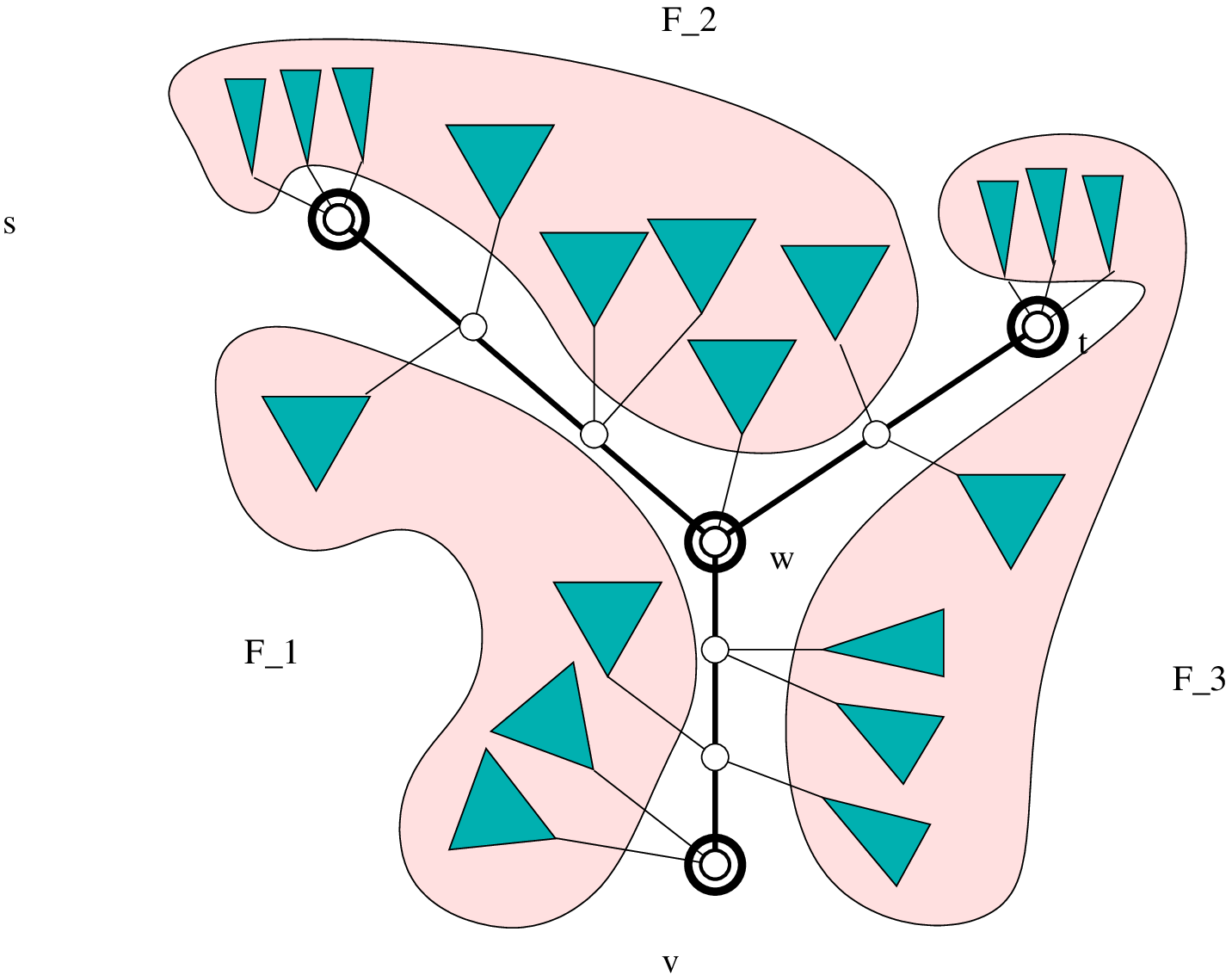}}
\captionn{\label{dectriple}On this example $w_1=321$, $w_2=12$,$w_3=2$, $l=1,r=3$, $S_1=4$, $S_2=8$, $S_3=7$.}
\end{figure}

To prove Proposition \ref{cvradon}, we now evaluate $`P(\Dec(n)=x)/`P(\widetilde{\Dec}(n)=x)$ for any $x\in A_{n,M}$. The number of ternary trees from $\Tter_{3n+1}$ satisfying
\[\Dec(n)
=(w_1,w_2,w_3,|w_1|,|w_2|,|w_3|,l,r)\] for some prescribed words $w_1,w_2,w_3$ and $(l,r)\in I_3$ is equal to the number of 3-tuples of forests as drawn on Figure \ref{dectriple}. The first forest $F_1$ has
$S_1(w_1,w_2,w_3,l,r) = N_1(w_1)+N_1(w_2)+l-1$ roots and since $w(ns)$ is the $ns+1$th internal nodes (not counted in $F_1$) and since the branch $\cro{\varnothing,w(ns)}$ contains $|w_1|+|w_2|+2$ internal nodes, $F_1$ has $n_1(w_1,w_2,w_3,l,r)= {ns}-|w_1|-|w_2|-1$ internal nodes (and then $3n_1+S_1$ nodes). The second forest $F_2$ has $S_2(w_1,w_2,w_3,l,r) = 3+N_2(w_2)+N_1(w_3)+(r-l-1)$ roots (the 3 comes from the fact that $w(ns)$ is an internal node), and $n_2(w_1,w_2,w_3,l,r)= {nt}-{ns}-|w_3|-1$ internal nodes. Finally the third forest $F_3$ has 
$S_3(w_1,w_2,w_3,l,r) = 3+ N_2(w_3)+N_2(w_1)+3-r$ roots and  $n_3(w_1,w_2,w_3,l,r)= n-nt-1$ internal nodes.
\par
\noindent Before going further, we recall that under $`P^\ter$ all trees in $\Tter_{3n+1}$ have the same weight $3^{-n}(2/3)^{2n+1}$ since they have $n$ internals nodes and $2n+1$ leaves. Let $F_k=(T(1),\dots,T(k))$ be a forest composed with $k$ independent GW trees with distribution $`P^\ter$, and let $|F_k|=\sum_{i=1}^k |T(i)|$ be the total number of nodes in $F_k$. By the rotation/conjugation principle,
\[`P^\ter(|F_k|=m)= \frac{k}{m}q({m,k})\] where $q({m,k})=`P(Z_m=-k)$ where $Z:=(Z_i)_{i\geq 0}$ is a random walk starting from 0, whose increment value are $-1$ or 2 with respective probability  $2/3$ and $1/3$.
\begin{lem}For any  words $w_1,w_2,w_3$ on the alphabet $\Sigma_3$ and $(l,r)\in I_3$, we have
  \be
`P^{\ter}_{3n+1}\l((W_1,W_2,W_3,L,R)=(w_1,w_2,w_3,l,r)\r)&=&
\frac{`P^\ter\l(|F^{i}_{S_i}|=3n_i+S_i, i\in\{1,2,3\}\r)}{3^{|w_1|+|w_2|+|w_3|+3}`P^{\ter}(\Tter_{3n+1})}\\
&=&
\frac{\prod_{i=1}^3 \frac{S_i}{3n_i+1}q_{3n_i+S_i,S_i} }{3^{|w_1|+|w_2|+|w_3|+3}`P^{\ter}(\Tter_{3n+1})}
\ee
where the $F^{i}$ are independent GW forests with respective number of roots the $S_i:=S_i(w_1,w_2,w_3,l,r)$'s, and  $n_i=n_i(w_1,w_2,w_3,l,r)$ for any $i\in\{1,2,3\}$.
\end{lem}
Notice that there is a hidden condition here since $(L,R)$  are well defined only when $u(ns)$ is not an ancestor of $u(nt)$ (which happens with probability going to 0).
\begin{note}\rm 
Notice that if $|w_i|=h_i$ for every $i$, for any $l,r\in\{1,2,3\}$, then 
\[`P^{\ter}_{3n+1}\l(\Dec(n)=(w_1,w_2,w_3,h_1,h_2,h_3,l,r)\r)=`P^{\ter}_{3n+1}\l((W_1,W_2,W_3,L,R)=(w_1,w_2,w_3,l,r)\r).\]
\end{note}
\proof This is just a counting argument, together with the remark that all the trees in $\Tter_{3n+1}$ have the same weight. The term $(1/3)^{|w_1|+|w_2|+|w_3|+3}$ comes from the $|w_1|+|w_2|+|w_3|+3$ internal nodes on the branches $\cro{\varnothing,w(ns)}\cup\cro{\varnothing,w(nt)}$. \hfill$\Box$
\medskip

We now evaluate $`P^{\ter}_{3n+1}\l(\widetilde{\Dec}(n)=(w_1,w_2,w_3,h_1,h_2,h_3,l,r)\r)$
for $(w_1,w_2,w_3)\in \Sigma_3^{h_1}\times\Sigma_3^{h_2}\times\Sigma_3^{h_3}$ and $(l,r)\in I_3$.  
The variable $\widetilde{\Dec}(n)$ is defined conditionally on $(H_1,H_2,H_3)$. We  have
\be
`P^{\ter}_{3n+1}\l((H_1,H_2,H_3)=(h_1,h_2,h_3)\r)&=&\sum
\frac{`P^\ter\l(|F^{i}_{S'_i}|=3n'_i+S'_i, i\in\{1,2,3\}\r)}{3^{|w'_1|+|w'_2|+|w'_3|+3}`P^{\ter}({\cal T}^\ter_{3n+1})}\\
&=&\sum\frac{\prod_{i=1}^3\frac{S'_i}{3n'_i+S'_i}q({3n'_i+S'_i,S'_i})}{3^{|w'_1|+|w'_2|+|w'_3|+3}`P^{\ter}({\cal T}^\ter_{3n+1})}
\ee
where $S_i':=S_i(w'_1,w'_2,w'_3,l',r')$'s, $n'_i=n_i(w'_1,w'_2,w'_3,l',r')$ and where the sum is taken on $(w'_1,w'_2,w'_3)\in \Sigma_3^{h_1}\times\Sigma_3^{h_2}\times\Sigma_3^{h_3}$ and $(l',r')\in I_3$. The term $3^{-|w'_1|-|w'_2|-|w'_3|-3}$ comes from the internal nodes of the branch  $\cro{\varnothing,w(ns)}\cup\cro{\varnothing,w(nt)}$. 
In other words
\begin{equation}
`P^{\ter}_{3n+1}\l((H_1,H_2,H_3)=(h_1,h_2,h_3)\r)=\frac{
`E\l(\prod_{i=1}^3\frac{{\bf S}_i}{3{\bf n}_i+{\bf S_i}}q(3{\bf n}_i+{\bf S}_i,{\bf S}_i)\r)}{3^2`P^{\ter}({\cal T}^\ter_{3n+1})}
\end{equation}
where ${\bf S}_i$ and ${\bf n}_i$ are the r.v. $S_i$ and $n_i$ when the $w_i$ are words with $h_i$ i.i.d. letters, uniform in $\Sigma_3$ and $({\bf l},{\bf r})$ is uniform in $I_3$.
Finally, by conditioning on the $H_i$'s, we get 
\[`P^{\ter}_{3n+1}\l(\widetilde{\Dec}(n)=(w_1,w_2,w_3,h_1,h_2,h_3,l,r)\r)=
\frac{`P^{\ter}_{3n+1}\l((H_1,H_2,H_3)=(h_1,h_2,h_3)\r)}{3^{|w_1|+|w_2|+|w_3|+1}}\]
and
\begin{equation}\label{rapp}
\frac{`P^{\ter}_{3n+1}\l(\Dec(n)=(w_1,w_2,w_3,h_1,h_2,h_3,l,r)\r)}
{`P^{\ter}_{3n+1}\l(\widetilde{\Dec}(n)=(w_1,w_2,w_3,h_1,h_2,h_3,l,r)\r)}=
\frac{\prod_{i=1}^3\frac{S_i}{3n_i+S_i}q({3n_i+S_i,S_i})}
{`E\l(\prod_{i=1}^3\frac{{\bf S}_i}{3{\bf n}_i+{\bf S}_i}q(3{\bf n}_i+{\bf S}_i,{\bf S}_i)\r)}.
\end{equation}
This quotient may be uniformly approached for $(w_1,w_2,w_3,h_1, h_2,h_3,l,r)\in A_{n,M}$ thanks to a central local limit theorem applied to the random walk $Z$:
\[\sup_{l\in -n+3\mathbb{N}}\l|\frac{\sqrt{n}}{3}\,`P(Z_n=l)-\frac{1}{\sqrt{4\pi}}\exp\l(-\frac{l^2}{4n}\r)\r|\xrightarrow[~~n~~]{} 0,\]
since the increment of $Z$ are centered and have variance 2.
This gives easily an equivalent for the numerator of \eref{rapp} (since
$q({m,k})=`P(Z_m=-k)$). For the denominator, split the expectation with respect to $(w'_1,w'_2,w'_3)$ belonging to $J_{h_1}\times J_{h_2}\times J_{h_3}$ or not. The first case occurs with probability close to 1, and the local central limit theorem provides the same asymptotic that the numerator. The second case provides an asymptotic with a smaller order (notice that the fact that $0\leq N_1(w)\leq 2|w|$ simplifies the use of the central local limit theorem) and we
get for any $`e>0$,
\[\l|\frac{`P^{\ter}_{3n+1}\l(\Dec(n)=(w_1,w_2,w_3,h_1,h_2,h_3,l,r)\r)}
{`P^{\ter}_{3n+1}\l(\widetilde{\Dec}(n)=(w_1,w_2,w_3,h_1,h_2,h_3,l,r)\r)}-1\r|\leq `e\]
on $A_{n,M}$ for  $n$ large enough.  $~\hfill\Box$

\section*{Acknowledgments} 
We would like to thank L. Devroye, D. Renault, N. Bonichon and O. Bernardi who pointed out several references.

\end{document}